\documentclass[11pt]{article}%
\usepackage[english]{babel}
\usepackage[utf8]{inputenc}
\usepackage[T1]{fontenc}
\usepackage{latexsym}
\usepackage[colorlinks=true, bookmarksnumbered=true,
bookmarksopen=true,bookmarksopenlevel=3, pdfstartview=FitH, linkcolor=cyan,
pdfmenubar=true,pdftoolbar=true, bookmarks=true,citecolor=cyan,
urlcolor=magenta,filecolor=cyan,menucolor=black,plainpages=false,pdfpagelabels]%
{hyperref}
\usepackage[lmargin=2cm,tmargin=2cm,bmargin=2cm,rmargin=2cm]{geometry}
\usepackage{amsbsy,amsmath,amsfonts,amssymb,amsthm}
\usepackage{graphicx}
\usepackage{graphics,xcolor}
\usepackage{amsmath}
\usepackage{amsfonts}
\usepackage{amssymb}%
\newtheorem{theo}{Theorem}[section]

\newtheorem{Remark}[theo]{Remark}
\newtheorem{lem}[theo]{Lemma}

\newtheorem{Proposition}[theo]{Proposition}
\numberwithin{equation}{section}

\newcommand\reallywidehat[1]{\savestack{\tmpbox}{\stretchto{  \scaleto{    \scalerel*[\widthof{\ensuremath{#1}}]{\kern-.6pt\bigwedge\kern-.6pt}    {\rule[-\textheight/2]{1ex}{\textheight}}  }{\textheight}}{0.5ex}}\stackon[1pt]{#1}{\tmpbox}}

\begin{document}

\title{Long-time solvability for the 2D inviscid Boussinesq equations \\with borderline regularity and dispersive effects}
\author{{V. Angulo-Castillo$^{1}$}\\{\small National University of Colombia, Campus Orinoquia, Department of
Mathematics,}\\{\small Kil\'{o}metro 9 v\'{\i}a a Ca\~{n}o Lim\'{o}n, Arauca, Colombia}\\{}\\{L. C. F. Ferreira$^{2}$}{\thanks{{Corresponding author. }\newline{E-mail
adresses: vlcastillo@unal.edu.co (V. Angulo-Castillo), lcff@ime.unicamp.br
(L.C.F. Ferreira), kosloff@ime.unicamp.br (L. Kosloff).} \newline{L.C.F.
Ferreira was supported by FAPESP (grant: 2020/05618-6) and CNPq (grant:
308799/2019-4), Brazil. L. Kosloff was supported by FAPESP (grant:
2016/15985-0), Brazil.}}} \ \ and \ \ {L. Kosloff}$^{3}$ \\{\small University of Campinas, IMECC-Department of Mathematics} \\{\small Rua S\'{e}rgio Buarque de Holanda, 651, CEP 13083-859, Campinas-SP,
Brazil}}
\date{}
\maketitle

\begin{abstract}
We are concerned with the long-time solvability for 2D inviscid Boussinesq
equations for a larger class of initial data which covers the case of
borderline regularity. First we show the local solvability in Besov spaces
uniformly with respect to a parameter $\kappa$ associated with the
stratification of the fluid. Afterwards, employing a blow-up criterion and
Strichartz-type estimates, the long-time solvability is obtained for large
$\kappa$ regardless of the size of initial data.

\bigskip{} \noindent\textbf{AMS MSC:} 35Q35; 76B03; 76U05; 35A01; 46E35

\medskip{} \noindent\textbf{Key:} Boussinesq equations; Convection problem;
Long-time solvability; Dispersive effects; Besov spaces; Borderline regularity

\end{abstract}

\section{Introduction}

We begin by considering the two-dimensional (2D) inviscid Boussinesq system%

\begin{equation}
\left\{
\begin{split}
&  \partial_{t}u+(u\cdot\nabla)u+\nabla p=\kappa\theta e_{2},\\
&  \partial_{t}\theta+(u\cdot\nabla)\theta=0,\\
&  \text{div}\,u=0,\\
&  u(x,0)=u_{0}(x),\;\theta(x,0)=\theta_{0}(x),
\end{split}
\right.  \label{Bouss0}%
\end{equation}
where $(x,t)\in\mathbb{R}^{2}\times(0,\infty)$, $u$ is the fluid velocity,
$\theta$ denotes the temperature (or the density in geophysical flows), $p$
stands for the pressure, $\kappa> 0$ is a gravitational constant which will be
associated with the stratification of the fluid (as described below), and the
vector $e_{2}=(0,1)$ indicates the positive vertical direction.

The 2D Boussinesq equations arise as a model in lower dimensions for the 3D
hydrodynamics equations by approximating the exact density of the fluid by a
constant representative value \cite{Salmon}. In particular, these equations
serve to model large scale atmospheric and oceanic flows that are responsible
for cold fronts and the jet stream (see \cite{Majda}, \cite{Pedlosky}).

Applying the \textquotedblleft curl\textquotedblright\ to the first equation
in (\ref{Bouss0}), and recalling that the vorticity $\omega= \text{curl}%
(u)=\nabla^{\perp}\cdot u$, we arrive at the equivalent vorticity formulation%

\begin{equation}
\left\{
\begin{split}
&  \partial_{t} \omega+(u\cdot\nabla) \omega= \kappa\partial_{1} \theta,\\
&  \partial_{t}\theta+(u\cdot\nabla)\theta=0,\\
&  u= \nabla^{\perp} (-\Delta)^{-1} \omega, \quad\; \nabla^{\perp}%
=(-\partial_{2}, \partial_{1}),\\
&  \omega(x,0)= \omega_{0}(x),\; \theta(x,0)=\theta_{0}(x),\quad\text{in
}\mathbb{R}^{2}.
\end{split}
\right.  \label{Bouss01}%
\end{equation}

Given that the atmosphere is physically observed to be mainly stable around
the hydrostatic balance between the pressure gradient and gravitational
effects \cite{Pedlosky}, we wish to consider the initial temperature close to
a physically nontrivial, stably stratified, stationary solution; namely, the
hydrostatic balance $\theta_{0}(x)=\rho_{0}(x)-\kappa x_{2}$, and then look
for the solution of (\ref{Bouss01}) with the temperature in the form
$\theta(t,x)=\rho(t,x)-\kappa x_{2}\,$. These modifications lead us to work
with the following new system (see \cite{Takada2021} for more details)%

\begin{equation}
\left\{
\begin{split}
&  \partial_{t} \omega+(u\cdot\nabla) \omega= \kappa\partial_{1} \rho,\\
&  \partial_{t}\rho+(u\cdot\nabla)\rho=\kappa u_{2},\\
&  u= \nabla^{\perp} (-\Delta)^{-1} \omega,\\
&  \omega(x,0)= \omega_{0}(x),\; \rho(x,0)=\rho_{0}(x),\quad\text{in
}\mathbb{R}^{2}.
\end{split}
\right.  \label{Bouss1}%
\end{equation}

The parameter $\kappa$ can then be interpreted as a gravitational constant
with $\kappa= N^{2}$, where $N>0$ is the buoyancy or Brunt-V\"{a}is\"{a}l\"{a}
frequency, representing the strength of stable stratification. System
(\ref{Bouss1}) thus exhibits a dispersive nature due to the stable
stratification terms, as will be developed below, and can be seen as a
Rayleigh-B\'{e}nard convection model where hot fluid sits on top of cooler fluid.

From the mathematical point of view, the inviscid 2D Boussinesq equations
(\ref{Bouss0}) are also important because they retain essential structural
features of the 3D Euler equations which derive from the vortex stretch
mechanism; in fact, they are equivalent to the 3D axisymmetric Euler equations
away from the symmetry axis. Additionally, this connection has prompted the
analysis of this system in the presence of nonlocal dissipative mechanisms for
the velocity or the temperature (or both), represented by a fractional
Laplacian operator. Global regularity has been shown in various scenarios of
these viscous systems, but there are still a number of open problems; see for
instance \cite{CCL}, \cite{SW}, \cite{WC17}, and references therein. In
particular, the challenging problem of stabilization around the hydrostatic
balance with partial viscosity or partial dissipation has been shown to be
feasible by specifically exploiting the wave structure provided by the
coupling of the velocity with the temperature (cf. \cite{T-W-Z-Z}, \cite{Wu
Stab}).

On the other hand, the regularity problem for the inviscid system
(\ref{Bouss0}) is in general more difficult, as it still appears to be unknown
whether the classical application of the Beale-Kato-Majda (BKM) regularity
criterion, based on the control of $\| \omega\|_{L^{\infty}}$, is enough for
the global regularity (cf. \cite{WuRef}). Thus, most results for this system
are local in time, with various works dealing with the local well-posedness,
in subcritical Besov and Sobolev spaces, by using the embedding of
$H^{s}(\mathbb{R}^{2})$ into $L^{\infty}$, $s>1$. The borderline case $s=1$
has added complications since $\| \nabla u \|_{L^{\infty}}$ is not bounded by
$\| \omega\|_{H^{1}}$, but this can be circumvented in the critical Besov
space to show local well-posedness, where the embedding into $L^{\infty}$ does
hold (cf. \cite{Liu et al}).

Moreover, system (\ref{Bouss0}) has been shown to be ill-posed in a borderline
or critical regularity setting, where the notion of criticality is defined as
the lower threshold where local well-posedness of strong solutions holds (cf.
\cite{Bourgain-Li},\cite{ElgMasm}). In particular, the approach in
\cite{ElgMasm} hinges on showing that ``the non-linearity does not serve as a
stabilizing mechanism" \cite{ElgMasm}, and is based on the condition that the
equation be locally well-posed in the critical Besov space which imbeds in
$L^{\infty}$, along with the application of a nontrivial commutator estimate;
so that the ill-posedness is linear and critical. Furthermore, it is shown
that system (\ref{Bouss1}) is ill-posed in the Yudovich class due to norm
inflation of the vorticity. Thus, these results provide the backdrop for the
ensuing proof of finite time blow-up for strong solutions of system
(\ref{Bouss0}) (cf. \cite{EJJ}).

However, there is a stark contrast for the inviscid system (\ref{Bouss1}) in
the presence of dispersive effects. The absence of dissipation presents the
possibility of instabilities, but in \cite{EW} Elgindi and Widmayer are able
to derive a sharp dispersive estimate and prove the long-time existence and
non-linear stability of (\ref{Bouss1}) about the stationary configuration, see
also \cite{Wan-CCM-2020} for related results. The approach consists in showing
that, in the case of $\kappa=1$, the dispersive estimate reveals the explicit
time decay rate of the linearized system, while not affecting the energy
estimates. The nonlinear stability then follows by improving the dependence of
the local time of existence on the size of the initial data. As this structure
is preserved at each level of the iterative scheme one can then obtain uniform
estimates for the approximation, as we show below. Deriving the corresponding
Strichartz estimates one can then control for the size of the initial data by
making $\kappa$ sufficiently large.

This dispersive estimate is also used for the analogous result in the case of
the dispersive inviscid surface quasigeostrophic (SQG) equation, where
dispersive effects are associated with strong rotation; indeed, both proofs
have the same structure. We refer the reader to \cite{AngFerr}, \cite{BaMaNi2}%
, \cite{Dutrifoy}, \cite{KLT}, and references therein, where the authors
showed that high speeds of rotation tend to smooth out 3D Navier-Stokes and
Euler flows.

In this direction, the long-time solvability of the system in the limit of
strong stratification is important to reveal the structure of solutions and
their main dynamics. For instance, in the case of the 3D Boussinesq system,
Widmayer \cite{Wid} shows that as the dispersive parameter grows to infinity,
the limiting system is a stratified system of 2D Euler equations with
stratified density.

This problem was analyzed for the SQG equations and the Boussinesq system
(\ref{Bouss1}) by Wan and Chen \cite{WC17} to show the long-time solvability
of strong solutions for large enough $\left\vert \kappa\right\vert $. The term
\textit{solvability} here refers to the pair of existence-uniqueness of
solutions in an appropriate sense. The approach in \cite{WC17} is based on
generalizing the dispersive estimate of \cite{EW} and deriving the
corresponding Strichartz estimates, to then show long-time solvability via a
blow-up criterion of BKM-type in Sobolev spaces $(\dot{H}^{s}\cap\dot{H}%
^{-1})\times H^{s+1}$, with $s>3$. Using an alternative argument, Takada
\cite{Takada2021} then improved this result with a weaker smoothness condition
on the initial data, as only belonging to $H^{s}$ ($s>2$), also showing the
asymptotics of solutions as $\kappa\rightarrow\infty$. Nevertheless, the
results of \cite{Takada2021} and \cite{WC17} do not reach the case $s=2$ that
appears as a borderline value for the long-time solvability of (\ref{Bouss1})
in Sobolev spaces $H_{p}^{s}$, and Besov spaces $B_{p,q}^{s}$, with $p=2$.

Thus, bearing in mind the lack of control on the vorticity in the context of
local well-posedness in the borderline Besov space, we are motivated to show
how the solvability can be extended in this context to arbitrary times through
the strong linear dispersive stabilization. For this we exploit the
paraproduct estimates which derive from the embedding of $\dot{B}_{2,1}^{1}$
into $L^{\infty}$, following the theme in Vishik's result of long-time
uniqueness for the 2D Euler equations in the borderline regularity case
\cite{Vishik99}. In particular, since we need to show local estimates which
are uniform in $\kappa$, in contrast to \cite{Liu et al}, we must develop
specific commutator estimates.

Then, in view of the identifications $B_{2,2}^{s}\equiv H^{s}$ and $\dot
{B}_{2,2}^{s}\equiv\dot{H}^{s}$ for $p=2$ and $q=2$ (see \cite[Theorem
6.4.4]{Berg-Lofstrom}), we are able to recover the results in \cite{WC17} and
also cover the case $s=2$, with similar power laws for the time decay rate.
Our main results read as follows:

\begin{theo}
\label{Main-theorem} Let $s$ and $q$ be real numbers such that $s>2$ with
$1\leq q\leq\infty$ or $s=2$ with $q=1$.

\begin{enumerate}
\item[$(i)$] (Local uniform solvability) Let $\omega_{0}\in\dot{B}_{2,q}%
^{s-1}(\mathbb{R}^{2})\cap\dot{H}^{-1}(\mathbb{R}^{2})$ and $\rho_{0}\in
B_{2,q}^{s}(\mathbb{R}^{2})$. There exists $T>0$ (depending of $\Vert
\omega_{0}\Vert_{\dot{B}_{2,q}^{s-1}\cap\dot{H}^{-1}}$ and $\Vert\rho_{0}%
\Vert_{B_{2,q}^{s}}$) such that (\ref{Bouss1}) has a unique solution
$(\omega,\rho)$ with $\omega\in C([0,T];\dot{B}_{2,q}^{s-1}(\mathbb{R}%
^{2})\cap\dot{H}^{-1}(\mathbb{R}^{2}))\cap C^{1}([0,T];\dot{B}_{2,q}%
^{s-2}(\mathbb{R}^{2})\cap\dot{H}^{-1}(\mathbb{R}^{2}))$ and $\rho\in
C([0,T];B_{2,q}^{s}(\mathbb{R}^{2}))\cap C^{1}([0,T];B_{2,q}^{s-1}%
(\mathbb{R}^{2}))$, for all $\kappa\in\mathbb{R}.$

\item[$(ii)$] (Long-time solvability) Let $T\in(0,\infty)$, $\omega_{0}\in
\dot{B}_{2,q}^{s}(\mathbb{R}^{2})\cap\dot{H}^{-1}(\mathbb{R}^{2})$, and
$\rho_{0}\in B_{2,q}^{s+1}(\mathbb{R}^{2})$. There exists $\kappa_{0}%
=\kappa_{0}(T,\Vert\omega_{0}\Vert_{\dot{B}_{2,q}^{s-1}\cap\dot{H}^{-1}}%
,\Vert\rho_{0}\Vert_{B_{2,q}^{s+1}})>0$ such that if $|\kappa|\geq\kappa_{0}$
then (\ref{Bouss1}) has a unique solution $(\omega,\rho)$ such that $\omega\in
C([0,T];\dot{B}_{2,q}^{s}(\mathbb{R}^{2})\cap\dot{H}^{-1}(\mathbb{R}^{2}))\cap
C^{1}([0,T];\dot{B}_{2,q}^{s-1}(\mathbb{R}^{2})\cap\dot{H}^{-1}(\mathbb{R}%
^{2}))$ and $\rho\in C([0,T];B_{2,q}^{s+1}(\mathbb{R}^{2}))\cap C^{1}%
([0,T];B_{2,q}^{s}(\mathbb{R}^{2}))$.
\end{enumerate}
\end{theo}

An analogous theorem for the SQG system was proved in our previous work
\cite{AngFerrKos}, and although the general argument is structured in a
similar vein, there are significant technical differences with the proof
provided here for Theorem \ref{Main-theorem}, which derive from the how the
stabilization due to the coupling of velocity and temperature works at
different levels.

In order to prove Theorem \ref{Main-theorem}, and show how the regularization
of the temperature feeds into further stabilization, we construct approximate
solutions $(\omega_{n},\rho_{n})_{n\in\mathbb{N}}$ via a Picard iteration
scheme, and show a priori estimates uniform with respect to the dispersive
parameter $\kappa$, so as to obtain a solution as the limit of $(\omega
_{n},\rho_{n})$ in the Besov spaces with the borderline regularity for
(\ref{Bouss1}) (see Section \ref{Localregularizationsolution} and the proof of
Theorem \ref{Main-theorem}).

To do this in the framework of borderline regularity, we employ an
intersection of spaces for the vorticity $\omega$. In fact, $\dot{H}^{-1}$ is
important to control the influence of $\kappa$ on the existence time $T$, and
then obtain a uniform time w.r.t $\kappa$, via a cancellation effect involving
the $\dot{H}^{-1}$-inner product (see, e.g., (\ref{cancelation})), while the
homogeneous Besov space $\dot{B}_{2,q}^{s-1}(\mathbb{R}^{2})$ provides the
necessary control on the regularity, particularly for the borderline case.
More precisely, we need to show that $(\omega_{n},\rho_{n})_{n\in\mathbb{N}}$
is bounded in $L^{\infty}(0,T;\dot{B}_{2,q}^{s-1}\cap\dot{H}^{-1})\times
L^{\infty}(0 ,T;B_{2,q}^{s})$ and Cauchy in $L^{\infty}(0,T;\dot{B}%
_{2,q}^{s-2}\cap\dot{H}^{-1})\times L^{\infty}(0 ,T;B_{2,q}^{s-1})$, both
uniformly w.r.t. $\kappa$. It is worth noting that this difficulty does not
appear in the context of inviscid SQG and Euler equations (with or without
dispersive effects) when analyzing local solvability and borderline regularity
in Besov spaces (see \cite{AngFerr,AngFerrKos,PakPark,Vishik99}).

We also note that to prove the uniform solvability in this context, it is
central to use commutator estimates in the framework of homogeneous Besov
spaces, which we present in a form which we could not find elsewhere in the
literature (see Section \ref{Section-3}). Subsequently, for large values of
$\left\vert \kappa\right\vert $ and $s\geq2$, we obtain the long-time
solvability by showing a blow-up criterion and handling globally the integral
$\int_{0}^{t}\Vert(\omega\pm\Lambda\rho)(\tau)\Vert_{\dot{B}_{\infty,1}^{0}%
}\,d\tau$\,, using Strichartz estimates as presented in \cite{KLT,WC17}.

Lastly, we remark that the stability problem for the more general setting
where the viscosity and thermal diffusivity are non-zero, which is also
physically important \cite{Gill}, \cite{Pedlosky}, has received considerable
attention. In particular, for the case of system \eqref{Bouss0} with positive
viscosity but zero thermal diffusivity, the global well-posedness and inviscid
limit have been shown through refined energy methods, see \cite{Chae},
\cite{H-K}, \cite{D-P}, \cite{LCF-EVR}, whereas the stabilization around the
hydrostatic balance has been recently studied in the case of bounded or strip
domains with Dirichlet or periodic boundary conditions, see \cite{CCL},
\cite{D-W-Z-Z}, \cite{T-W-Z-Z}. For the latter problem, the viscous term
introduces a wave structure in the linearized system for which explicit decay
rates can be shown, but so that the convergence to the full nonlinear
equations is considerably more involved. Given that the energy estimates for
this system has a parallel structure to the inviscid system, it is possible
that the iterative scheme could be useful for controlling the decay rate and
obtaining the asymptotic limit of large dispersive forcing.

The plan of the manuscript is as follows. In Section 2 we present some
preliminaries about Besov spaces and Strichartz estimates, among others.
Section 3 is devoted to the commutator estimates. In Section 4, we analyze the
approximation scheme $(\omega_{n},\rho_{n})_{n\in\mathbb{N}}$\ and obtain the
local-in-time solvability of (\ref{Bouss1}) uniformly with respect to the
parameter $\kappa$. The proof of Theorem \ref{Main-theorem} (ii) is carried
out in Section 5.

\section{Preliminaries}

The purpose of this section is to provide some basic definitions and
properties about Besov spaces as well as some estimates useful for our ends,
such as product, embeddings, Strichartz estimates, among others. We refer the
reader to the book \cite{Berg-Lofstrom} for more details on Besov spaces and
their properties.

First, we denote the Schwartz space on $\mathbb{R}^{2}$ by $\mathcal{S}%
(\mathbb{R}^{2})$ and its dual by $\mathcal{S}^{\prime}(\mathbb{R}^{2})$
(tempered distributions). For $f\in\mathcal{S}^{\prime}(\mathbb{R}^{2})$,
$\widehat{f}$ stands for the Fourier transform of $f$. Select a radial
function $\psi_{0}\in\mathcal{S}(\mathbb{R}^{2})$ satisfying $0\leq
\widehat{\psi}_{0}(\xi)\leq1$, $\mbox{supp}\,\widehat{\psi}_{0}\subset\left\{
\xi\in\mathbb{R}^{2}\colon\frac{5}{8}\leq|\xi|\leq\frac{7}{4}\right\}  $ and
\[
\sum_{j\in\mathbb{Z}}\widehat{\psi}_{j}(\xi)=1\quad\text{for all}\quad\xi
\in\mathbb{R}^{2}\setminus\{0\},
\]
where $\psi_{j}(x):=2^{2j}\psi_{0}(2^{j}x)$. For each $k\in\mathbb{Z}$, we
consider $S_{k},\dot{S}_{k}\in\mathcal{S}$ defined in Fourier variables as
\[
\widehat{S_{k}}(\xi)=1-\sum_{j\geq k+1}\widehat{\psi}_{j}(\xi)\quad
\text{and}\quad\widehat{\dot{S}_{k}}(\xi)=\sum_{j\leq k}\widehat{\psi}_{j}%
(\xi).
\]
We observe that
\[
\mbox{supp}\,\widehat{\psi}_{j}\cap\mbox{supp}\,\widehat{\psi}_{j^{\prime}%
}=\emptyset\quad\text{if}\quad|j-j^{\prime}|\geq2.
\]
For $f\in\mathcal{S}^{\prime}(\mathbb{R}^{2})$ and $j\in\mathbb{Z}$, the
Littlewood-Paley operator $\Delta_{j}$ is the convolution $\Delta_{j}%
f:=\psi_{j}\ast f$ which works as a filter on the support of $\psi_{j}.$ We
also consider the family of operators $\{\tilde{\Delta}_{k}\}_{k\in
\{0\}\cup\mathbb{N}}$ defined as $\tilde{\Delta}_{0}=S_{0}\ast f$ and
$\tilde{\Delta}_{k}=\Delta_{k}$ for every integer $k\geq1.$

Let $\mathcal{P}=\mathcal{P}(\mathbb{R}^{2})$ denote the set of polynomials
and consider $s\in\mathbb{R}$ and $p,q\in\lbrack1,\infty].$ The homogeneous
Besov space $\dot{B}_{p,q}^{s}(\mathbb{R}^{2})$ is the set of all
$f\in\mathcal{S}^{\prime}(\mathbb{R}^{2})/\mathcal{P}$ such that
\[
\Vert f\Vert_{\dot{B}_{p,q}^{s}}:=\left\Vert \left\{  2^{sj}\Vert\Delta
_{j}f\Vert_{L^{p}}\right\}  \right\Vert _{l^{q}(\mathbb{Z})}<\infty.
\]
The nonhomogeneous version of $\dot{B}_{p,q}^{s}(\mathbb{R}^{2})$, namely the
Besov space $B_{p,q}^{s}=B_{p,q}^{s}(\mathbb{R}^{2}),$ is the space of all
$f\in\mathcal{S}^{\prime}(\mathbb{R}^{n})$ such that the norm $\left\Vert
f\right\Vert _{B_{p,q}^{s}}<\infty,$ where
\[
\left\Vert f\right\Vert _{B_{p,q}^{s}}=\left\{
\begin{array}
[c]{rll}%
\left(  \sum_{k=0}^{\infty}2^{ksq}\left\Vert \tilde{\Delta}_{k}f\right\Vert
_{p}^{q}\right)  ^{\frac{1}{q}}, & \text{ if }q<\infty, & \\
\sup_{k\in\{0\}\cup\mathbb{N}}\{2^{ks}\left\Vert \tilde{\Delta}_{k}%
f\right\Vert _{p}\}, & \text{ if \ }q=\infty. &
\end{array}
\right.  \label{norm-Besov-1}%
\]
The pairs $(\dot{B}_{p,q}^{s},\Vert\cdot\Vert_{\dot{B}_{p,q}^{s}})$ and
$(B_{p,q}^{s},\Vert\cdot\Vert_{B_{p,q}^{s}})$ are Banach spaces. Also, for
$s\in\mathbb{R}$ and $p,q\in\lbrack1,\infty]$, it follows that%
\[
\Vert f\Vert_{B_{p,q}^{s}}\leq C\left(  \Vert f\Vert_{L^{p}}+\Vert
f\Vert_{\dot{B}_{p,q}^{s}}\right)  .\label{Besov-equiv}%
\]
For $s>0,$ we have the equivalence of norms%
\[
\Vert f\Vert_{B_{p,q}^{s}}\sim\Vert f\Vert_{L^{p}}+\Vert f\Vert_{\dot{B}%
_{p,q}^{s}}.\label{Besov-equiv-2}%
\]
In the case $s=0,$ we recall the inclusion $B_{p,1}^{0}\hookrightarrow L^{p},$
for all $p\in\lbrack1,\infty].$

\begin{lem}
[Bernstein's Lemma]\label{inequalitybernstein} Let $1\leq p\leq\infty$ and
$f\in L^{p}$ be such that $\mbox{supp}\ \widehat{f}\subset\{\xi\in
\mathbb{R}^{2}:2^{j-2}\leq|\xi|\leq2^{j}\}$. Then, we have the estimates
\[
C^{-1}2^{jk}\Vert f\Vert_{L^{p}}\leq\Vert D^{k}f\Vert_{L^{p}}\leq C2^{jk}\Vert
f\Vert_{L^{p}},
\]
where $C=C(k)$ is a positive constant.
\end{lem}

\begin{Remark}
\label{remark1} Using the above lemma, one can prove the equivalence
\[
\Vert f\Vert_{{\dot{B}}_{p,q}^{s+k}}\sim\Vert D^{k}f\Vert_{{\dot{B}}_{p,q}%
^{s}}.
\]
Moreover, considering $1\leq p,q\leq\infty$ and $s\geq n/p,$ with $q=1$ if
$s=n/p,$ we have that (see, e.g., \cite[Sections 6.5 and 6.8]{Berg-Lofstrom})
\begin{equation}
\label{nonhom-ineq-1}\Vert f\Vert_{L^{\infty}}\leq C\Vert f\Vert_{B_{p,q}^{s}%
}.
\end{equation}
Then,
\begin{equation}
\label{nonhom-ineq-2}\Vert\nabla f\Vert_{L^{\infty}}\leq C\Vert\nabla
f\Vert_{B_{p,q}^{s-1}}\leq C\Vert f\Vert_{B_{p,q}^{s}},
\end{equation}
where $1\leq p,q\leq\infty$ and $s\geq n/p+1$ with $q=1$ in the case $s=n/p+1$.

The uniform estimates we develop rely on inequalities (\ref{nonhom-ineq-1})
and (\ref{nonhom-ineq-2}) in their respective regularity ranges $s\geq n/p$
and $s\geq n/p+1$. For the case $q=1$ and the corresponding $s$, these
inequalities hold in homogeneous Besov spaces. However, as our results are
aimed for more general $q \geq1$ and $s$, and deal with the full inviscid
case, we restrict to the case of nonhomogeneous spaces, where they are
guaranteed to be valid, and leave open the question of working with
homogeneous Besov spaces for the component $\rho$.
\end{Remark}

Some Leibniz-type rules in Besov spaces are the subject of the lemma below
(see \cite{Chae2004}).

\begin{lem}
\label{inequalityholder} Let $s>0$, $1\leq p_{1},p_{2}\leq\infty$, $1\leq
r_{1},r_{2}\leq\infty$ and $1\leq p,q\leq\infty$ be such that $\frac{1}%
{p}=\frac{1}{p_{1}}+\frac{1}{p_{2}}=\frac{1}{r_{1}}+\frac{1}{r_{2}}.$ Then, we
have the estimates
\[%
\begin{split}
\Vert fg\Vert_{{\dot{B}}_{p,q}^{s}}  &  \leq C(\Vert g\Vert_{L^{p_{2}}}\Vert
f\Vert_{{\dot{B}}_{p_{1},q}^{s}}+\Vert f\Vert_{L^{r_{2}}}\Vert g\Vert
_{{\dot{B}}_{r_{1},q}^{s}}),\\
\Vert fg\Vert_{B_{p,q}^{s}}  &  \leq C(\Vert g\Vert_{L^{p_{2}}}\Vert
f\Vert_{B_{p_{1},q}^{s}}+\Vert f\Vert_{L^{r_{2}}}\Vert g\Vert_{B_{r_{1},q}%
^{s}}),
\end{split}
\]
where $C>0$ is a universal constant.
\end{lem}

We will employ the Strichartz estimates of \cite{WC17}, linked to the
dispersive term $\kappa(\partial_{1}\rho+\Lambda u_{2})$ obtained from
(\ref{Bouss1}), which will allow us to obtain long-time solvability for
(\ref{Bouss1}). In particular, we will use the following results found in
\cite{AngFerrKos}, \cite{KLT} and \cite{WC17}.

\begin{lem}
\label{LemmaStrichartz} Let $\kappa\in\mathbb{R}$, $4\leq\gamma\leq\infty$ and
$2\leq r\leq\infty$ be such that
\begin{equation}
\frac{1}{\gamma}+\frac{1}{2r}\leq\frac{1}{4}. \label{condition_to_gamma_and_r}%
\end{equation}
Then, there holds
\[
\left\Vert \mathcal{G}_{\pm}(\kappa t)f\right\Vert _{L^{\gamma}(0,\infty
;L^{r})}\leq C|\kappa|^{-\frac{1}{\gamma}}\Vert f\Vert_{L^{2}},
\]
for all $f\in L^{2}(\mathbb{R}^{2})$, where
\[
\mathcal{G}_{\pm}(t)f(x):=\int_{\mathbb{R}^{2}}e^{i\xi\cdot x\pm it\frac
{\xi_{1}}{|\xi|}}\widehat{\phi}(\xi)\widehat{f}(\xi)\ d\xi
\]
and $\widehat{\phi}$ is a compactly supported smooth function in
$\mathbb{R}^{2}$.
\end{lem}

\begin{lem}
\label{theorem-besov-strichartz} Let $s,t,\kappa\in\mathbb{R}$, $1\leq
q\leq\infty$, $4\leq\gamma\leq q,$ and $2\leq r\leq\infty$. Assume also
(\ref{condition_to_gamma_and_r}). Then,
\[
\left\Vert e^{\pm t\kappa\mathcal{R}_{1}}f\right\Vert _{L^{\gamma}%
(0,\infty;\dot{B}_{r,q}^{s})}\leq C|\kappa|^{-\frac{1}{\gamma}}\Vert
f\Vert_{\dot{B}_{2,q}^{s+1-\frac{2}{r}}}, \label{Strichartz_general}%
\]
for all $f\in\dot{B}_{2,q}^{s+1-\frac{2}{r}}(\mathbb{R}^{2})$.
\end{lem}

\section{Commutator estimates}

\label{Section-3}The present section is devoted to commutator estimates in
$\dot{B}_{p,q}^{s}$ and $B_{p,q}^{s}$ that will be useful to obtain
convergence of our approximate solutions. We state and prove some of them, as
we have not been able to locate them in the literature with the needed
hypotheses and conclusions for our purposes.

Recall the commutator operator
\[
\lbrack f\cdot\nabla,\Delta_{j}]g=f\cdot\nabla(\Delta_{j}g)-\Delta_{j}%
(f\cdot\nabla g).\label{aux-commut-1}%
\]
Using the H\"{o}lder inequality in a slightly different way than used in
\cite{Chae2004,Takada2008,WuXuYe2015}, it is possible to obtain the following
estimates for the commutator:

\begin{lem}
\label{2.4} Let $1<p<\infty$ and $1\leq q,p_{1},p_{2},r_{1},r_{2}\leq\infty$
be such that $\frac{1}{p}=\frac{1}{p_{1}}+\frac{1}{p_{2}}=\frac{1}{r_{1}%
}+\frac{1}{r_{2}}$.

\begin{enumerate}
\item[$(i)$] Let $s>0$, $f\in{\dot{B}}_{p_{1},q}^{s}(\mathbb{R}^{n})$ and
$g\in{\dot{B}}_{r_{1},q}^{s}(\mathbb{R}^{n})$. Assume further that $\nabla
f\in L^{r_{2}}(\mathbb{R}^{n})$\thinspace$,$ $\nabla\cdot f=0$ and $\nabla
g\in L^{p_{2}}(\mathbb{R}^{n})$. Then, we have the estimate
\[
\left(  \sum_{j\in\mathbb{Z}}2^{sjq}\Vert\lbrack f\cdot\nabla,\Delta
_{j}]g\Vert_{L^{p}}^{q}\right)  ^{1/q}\leq C\left(  \Vert\nabla f\Vert
_{L^{r_{2}}}\Vert g\Vert_{{\dot{B}}_{r_{1},q}^{s}}+\Vert\nabla g\Vert
_{L^{p_{2}}}\Vert f\Vert_{{\dot{B}}_{p_{1},q}^{s}}\right)  ,
\]
where $C>0$ is a universal constant.

\item[$(ii)$] Let $s>-1$, $f\in{\dot{B}}_{p_{1},q}^{s+1}(\mathbb{R}^{n})$ and
$g\in{\dot{B}}_{r_{1},q}^{s}(\mathbb{R}^{n})\cap L^{p_{2}}(\mathbb{R}^{n}).$
Assume further that $\nabla f\in L^{r_{2}}(\mathbb{R}^{n})$ and $\nabla\cdot
f=0$. Then, we have the estimate
\[
\left(  \sum_{j\in\mathbb{Z}}2^{sjq}\Vert\lbrack f\cdot\nabla,\Delta
_{j}]g\Vert_{L^{p}}^{q}\right)  ^{1/q}\leq C\left(  \Vert\nabla f\Vert
_{L^{r_{2}}}\Vert g\Vert_{{\dot{B}}_{r_{1},q}^{s}}+\Vert g\Vert_{L^{p_{2}}%
}\Vert f\Vert_{{\dot{B}}_{p_{1},q}^{s+1}}\right)  ,
\]
where $C>0$ is a universal constant.
\end{enumerate}
\end{lem}

\begin{Remark}
\label{remark2} With the help of Lemma 6.3 in \cite{WuXuYe2015} (see also
\cite{BahouriCheminDanchin}), the properties of the operator $\Lambda^{-1}$,
considering the same hypotheses of Lemma \ref{2.4} and using the same
arguments in the proof of the same result, we obtain the following commutator
estimate
\[
\left(  \sum_{j\in\mathbb{Z}}2^{sjq}\Vert\lbrack f\cdot\nabla,\Lambda
^{-1}\Delta_{j}]g\Vert_{L^{p}}^{q}\right)  ^{1/q}\leq C\left(  \Vert\nabla
f\Vert_{L^{r_{2}}}\Vert g\Vert_{{\dot{B}}_{r_{1},q}^{s-1}}+\Vert
g\Vert_{L^{p_{2}}}\Vert f\Vert_{{\dot{B}}_{p_{1},q}^{s}}\right)  .
\]

\end{Remark}

Recall that the Bony formula for the paraproduct of $f$ and $g$ is given by
\begin{equation}
fg=T_{f}g+T_{g}f+R(f,g), \label{aux-paraproduct-0}%
\end{equation}
where
\begin{equation}
T_{f}g:=\sum_{j\in\mathbb{Z}}S_{j-2}f\Delta_{j}g\quad\text{and}\quad
R(f,g):=\sum_{|j-j^{\prime}|\leq1}\Delta_{j}f\Delta_{j^{\prime}}g.
\label{aux-paraproduct-1}%
\end{equation}
In the sequel we state and prove the following commutator-type estimates:

\begin{lem}
\label{estimate_type_commutator} Let $1<p<\infty$ and $1\leq q,p_{1}%
,p_{2},r_{1},r_{2}\leq\infty$ be such that $\frac{1}{p}=\frac{1}{p_{1}}%
+\frac{1}{p_{2}}=\frac{1}{r_{1}}+\frac{1}{r_{2}}$.

\begin{enumerate}
\item[$(i)$] Let $s>0$, $f\in\dot{B}_{p_{1},q}^{s}(\mathbb{R}^{n})$ with
$\nabla f\in L^{r_{2}}(\mathbb{R}^{n})$ and $\nabla\cdot f=0$, and $g\in
\dot{B}_{r_{1},q}^{s}(\mathbb{R}^{n})$ with $\nabla g\in L^{p_{2}}%
(\mathbb{R}^{n})$. Then, there exists a universal constant $C>0$ such that
\[
\left(  \sum_{j\in\mathbb{Z}}2^{sjq}\left\Vert (\dot{S}_{j-2}f\cdot
\nabla)\Delta_{j}g-\Delta_{j}(f\cdot\nabla)g\right\Vert _{L^{p}}^{q}\right)
^{1/q}\leq C\left(  \Vert\nabla f\Vert_{L^{r_{2}}}\Vert g\Vert_{\dot{B}%
_{r_{1},q}^{s}}+\Vert\nabla g\Vert_{L^{p_{2}}}\Vert f\Vert_{\dot{B}_{p_{1}%
,q}^{s}}\right)  .
\]

\item[$(ii)$] Let $s>-1$, $f\in\dot{B}_{p_{1},q}^{s+1}(\mathbb{R}^{n})$ with
$\nabla f\in L^{r_{2}}(\mathbb{R}^{n})$ and $\nabla\cdot f=0$, and $g\in
\dot{B}_{r_{1},q}^{s}(\mathbb{R}^{n})\cap L^{p_{2}}(\mathbb{R}^{n})$. Then,
there exists a universal constant $C>0$ such that
\[
\left(  \sum_{j\in\mathbb{Z}}2^{sjq}\left\Vert (\dot{S}_{j-2}f\cdot
\nabla)\Delta_{j}g-\Delta_{j}(f\cdot\nabla)g\right\Vert _{L^{p}}^{q}\right)
^{1/q}\leq C\left(  \Vert\nabla f\Vert_{L^{r_{2}}}\Vert g\Vert_{\dot{B}%
_{r_{1},q}^{s}}+\Vert g\Vert_{L^{p_{2}}}\Vert f\Vert_{\dot{B}_{p_{1},q}^{s+1}%
}\right)  .
\]

\end{enumerate}
\end{lem}

\textbf{Proof.} The proof of part $(i)$, it follows from the calculations
obtained by Chae in \cite{Chae2004}. We show the part $(ii)$. We follow
closely the argument in \cite{Vishik99} (see also
\cite{Bahouri-Chemin94,Chemin92}). By Bony's paraproduct formula
(\ref{aux-paraproduct-0}), we can write
\[%
\begin{split}
(\dot{S}_{j-2}f\cdot\nabla)\Delta_{j}g-\Delta_{j}(f\cdot\nabla)g=  &
-\sum_{i=1}^{n}\Delta_{j}T_{\partial_{i}g_{k}}f_{i}+\sum_{i=1}^{n}[T_{f_{i}%
}\partial_{i},\Delta_{j}]g_{k}-\sum_{i=1}^{n}\Delta_{j}T_{f_{i}-S_{j-2}f_{i}%
}(\partial_{i}\Delta_{j}g_{k})\\
&  -\sum_{i=1}^{n}\left\{  \Delta_{j}R(f_{i},\partial_{i}g_{k})-R(S_{j-2}%
f_{i},\partial_{i}\Delta_{j}g_{k})\right\} \\
:=  &  I+II+III+IV.
\end{split}
\]
For $I$, in view of (\ref{aux-paraproduct-1}), it follows that
\[
I=-\sum_{j^{\prime}}\sum_{i=1}^{n}\Delta_{j}{\dot{S}_{j^{\prime}-2}%
(\partial_{i}g_{k})\Delta_{j^{\prime}}f_{i}}.
\]
We observe that $\text{supp}\mathcal{F}(\dot{S}_{j^{\prime}-2}(\partial
_{i}g_{k})\Delta_{j^{\prime}}f_{i})\subset\{\xi:2^{j^{\prime}-3}\leq|\xi
|\leq2^{j^{\prime}+1}\}$ and $\Delta_{j}{\dot{S}_{j^{\prime}-2}(\partial
_{i}g_{k})\Delta_{j^{\prime}}f_{i}}=0$, if $|j-j^{\prime}|\geq4$. Then,
\[
I=-\sum_{|j-j^{\prime}|\leq3}\sum_{i=1}^{n}\Delta_{j}{\dot{S}_{j^{\prime}%
-2}(\partial_{i}g_{k})\Delta_{j^{\prime}}f_{i}}.
\]
Using integration by parts, we arrive at
\[%
\begin{split}
I  &  =-\sum_{|j-j^{\prime}|\leq3}\sum_{i=1}^{n}2^{jn}\int_{\mathbb{R}^{n}%
}\psi_{0}(2^{j}(x-y))(\dot{S}_{j^{\prime}-2}\partial_{i}g_{k})(y)(\Delta
_{j^{\prime}}f_{i})(y)\,dy\\
&  =-\sum_{|j-j^{\prime}|\leq3}\sum_{i=1}^{n}2^{j}2^{jn}\int_{\mathbb{R}^{n}%
}\partial_{i}\psi_{0}(2^{j}(x-y))(\dot{S}_{j^{\prime}-2}g_{k})(y)(\Delta
_{j^{\prime}}f_{i})(y)\,dy\\
&  =-\sum_{|j-j^{\prime}|\leq3}\sum_{i=1}^{n}2^{j}\{(2^{jn}\partial_{i}%
\psi_{0}(2^{j}\cdot)\ast((\dot{S}_{j^{\prime}-2}g_{k})(\Delta_{j^{\prime}%
}f_{i}))\},
\end{split}
\]
which yields
\begin{equation}%
\begin{split}
\Vert I\Vert_{L^{p}}  &  \leq C\sum_{|j-j^{\prime}|\leq3}\sum_{i=1}^{n}%
2^{j}\Vert(\dot{S}_{j^{\prime}-2}g_{k})(\Delta_{j^{\prime}}f_{i})\Vert_{L^{p}%
}\\
&  \leq C\Vert g\Vert_{L^{p_{2}}}\sum_{|j-j^{\prime}|\leq3}2^{j}\Vert
\Delta_{j^{\prime}}f\Vert_{L^{p_{1}}}.
\end{split}
\label{estimate_I}%
\end{equation}
For estimate $II$, by an argument similar to the one above, we first note
that
\[
\lbrack\dot{S}_{j^{\prime}-2}f_{i},\Delta_{j}](\partial_{i}\Delta_{j^{\prime}%
}g_{k})=0,\text{ if }\ |j-j^{\prime}|\geq4.
\]
Then, using $\nabla\cdot\dot{S}_{j^{\prime}-2}f=0$ and integration by parts,
it holds that
\[%
\begin{split}
II  &  =\sum_{i=1}^{n}\sum_{|j-j^{\prime}|\leq3}\{(\dot{S}_{j^{\prime}-2}%
f_{i})\Delta_{j}(\partial_{i}\Delta_{j^{\prime}}g_{k})-\Delta_{j}(\dot
{S}_{j^{\prime}-2}f_{i})(\partial_{i}\Delta_{j^{\prime}}g_{k})\}\\
&  =\sum_{i=1}^{n}\sum_{|j-j^{\prime}|\leq3}2^{jn}\int_{\mathbb{R}^{n}}%
\psi_{0}(2^{j}(x-y))(\dot{S}_{j^{\prime}-2}f_{i}(x)-\dot{S}_{j^{\prime}%
-2}f_{i}(y))(\partial_{i}\Delta_{j^{\prime}}g_{k})(y)\,dy\\
&  =\sum_{i=1}^{n}\sum_{|j-j^{\prime}|\leq3}2^{j(n+1)}\int_{\mathbb{R}^{n}%
}\partial_{i}\psi_{0}(2^{j}(x-y))(\dot{S}_{j^{\prime}-2}f_{i}(x)-\dot
{S}_{j^{\prime}-2}f_{i}(y))(\Delta_{j^{\prime}}g_{k})(y)\,dy\\
&  =\sum_{i=1}^{n}\sum_{|j-j^{\prime}|\leq3}2^{j(n+1)}\int_{\mathbb{R}^{n}%
}\partial_{i}\psi_{0}(2^{j}(x-y))\int_{0}^{1}((x-y)\cdot\nabla)(\dot
{S}_{j^{\prime}-2}f_{i}(x+\tau(y-x)))\,d\tau\,(\Delta_{j^{\prime}}%
g_{k})(y)\,dy\\
&  =\sum_{i=1}^{n}\sum_{|j-j^{\prime}|\leq3}\int_{\mathbb{R}^{n}}\partial
_{i}\psi_{0}(z)\int_{0}^{1}(z\cdot\nabla)(\dot{S}_{j^{\prime}-2}f_{i}%
(x+\tau2^{-j}z))\,d\tau\,(\Delta_{j^{\prime}}g_{k})(x-2^{-j}z)\,dz.\\
&
\end{split}
\]
Therefore,
\[%
\begin{split}
|II|  &  \leq\sum_{i=1}^{n}\sum_{|j-j^{\prime}|\leq3}\int_{\mathbb{R}^{n}%
}|\partial_{i}\psi_{0}(z)|\int_{0}^{1}|z||\nabla(\dot{S}_{j^{\prime}-2}%
f_{i}(x+\tau2^{-j}z))|\,d\tau\,|(\Delta_{j^{\prime}}g_{k})(x-2^{-j}z)|\,dz\\
&  \leq C\Vert\nabla f\Vert_{L^{\infty}}\sum_{|j-j^{\prime}|\leq3}%
\int_{\mathbb{R}^{n}}|z||\nabla\psi_{0}(z)||(\Delta_{j^{\prime}}%
g_{k})(x-2^{-j}z)|\,dz,
\end{split}
\]
which leads us to
\begin{equation}%
\begin{split}
\Vert II\Vert_{L^{p}}  &  \leq C\Vert\nabla f\Vert_{L^{\infty}}\sum
_{|j-j^{\prime}|\leq3}\int_{\mathbb{R}^{n}}|z||\nabla\psi_{0}(z)|\Vert
(\Delta_{j^{\prime}}g_{k})(\cdot-2^{-j}z)\Vert_{L^{p}}\,dz\\
&  \leq C\Vert\nabla f\Vert_{L^{r_{1}}}\sum_{|j-j^{\prime}|\leq3}\Vert
\Delta_{j^{\prime}}g\Vert_{L^{r_{2}}}.
\end{split}
\label{estimate_II}%
\end{equation}
For $III$, we have that%
\[%
\begin{split}
III  &  =\sum_{i=1}^{n}\sum_{|j-j^{\prime}|\leq1}\dot{S}_{j^{\prime}-2}%
(f_{i}-\dot{S}_{j-2}f_{i})\partial_{i}\Delta_{j^{\prime}}\Delta_{j}g_{k}\\
&  =\sum_{i=1}^{n}\sum_{|j-j^{\prime}|\leq1}\dot{S}_{j^{\prime}-2}\left(
\sum_{m=j-1}^{j^{\prime}-1}\Delta_{m}f_{i}\right)  \partial_{i}\Delta
_{j}\Delta_{j^{\prime}}g_{k}\\
&  =\sum_{i=1}^{n}\sum_{|j-j^{\prime}|\leq1}\dot{S}_{j^{\prime}-2}\left(
\Delta_{j-1}f_{i}+\Delta_{j}f_{i}\right)  \partial_{i}\Delta_{j}%
\Delta_{j^{\prime}}g_{k}.
\end{split}
\]
Applying the $L^{p}$-norm, we arrive at
\begin{equation}%
\begin{split}
\Vert III\Vert_{L^{p}}  &  \leq\sum_{i=1}^{n}\sum_{|j-j^{\prime}|\leq1}%
(\Vert\Delta_{j-1}f_{i}\Vert_{L^{\infty}}+\Vert\Delta_{j}f_{i}\Vert
_{L^{\infty}})\Vert\partial_{i}\Delta_{j}\Delta_{j^{\prime}}g_{k}\Vert_{L^{p}%
}\\
&  \leq\sum_{i=1}^{n}\sum_{|j-j^{\prime}|\leq1}(2^{-j+1}\Vert\Delta
_{j-1}\nabla f_{i}\Vert_{L^{\infty}}+2^{-j}\Vert\Delta_{j}\nabla f_{i}%
\Vert_{L^{\infty}})2^{j}\Vert\Delta_{j}\Delta_{j^{\prime}}g_{k}\Vert_{L^{p}}\\
&  \leq C\Vert\nabla f\Vert_{L^{r_{1}}}\sum_{|j-j^{\prime}|\leq1}\Vert
\Delta_{j^{\prime}}g\Vert_{L^{r_{2}}}.
\end{split}
\label{estimate_III}%
\end{equation}
For the parcel $IV,$ we can decompose
\[%
\begin{split}
IV  &  =\sum_{i=1}^{n}\Delta_{j}\partial_{i}R(f_{i}-\dot{S}_{j-2}f_{i}%
,g_{k})+\sum_{i=1}^{n}\{\Delta_{j}R(\dot{S}_{j-2}f_{i},\partial_{i}%
g_{k})-R(\dot{S}_{j-2}f_{i},\Delta_{j}\partial_{i}g_{k})\}\\
&  =IV_{1}+IV_{2}.
\end{split}
\]
Since $\sum_{i=1}^{n}\partial_{i}\Delta_{j^{\prime}}(f_{i}-\dot{S}_{j-2}%
f_{i})=0 $, it follows that
\[%
\begin{split}
IV_{1}  &  =\sum_{i=1}^{n}\partial_{i}\Delta_{j}\left\{  \sum_{|j^{\prime
}-j^{\prime\prime}|\leq1}\Delta_{j^{\prime}}(f_{i}-\dot{S}_{j-2}f_{i}%
)\Delta_{j^{\prime\prime}}g_{k}\right\} \\
&  =\sum_{i=1}^{n}\Delta_{j}\left\{  \sum_{|j^{\prime}-j^{\prime\prime}|\leq
1}\sum_{j^{\prime}\geq j-3}(\Delta_{j^{\prime}}f_{i}-\dot{S}_{j-2}%
\Delta_{j^{\prime}}f_{i})\Delta_{j^{\prime\prime}}\partial_{i}g_{k}\right\} \\
&  =\sum_{i=1}^{n}\sum_{|j^{\prime}-j^{\prime\prime}|\leq1}\sum_{j^{\prime
}\geq j-3}2^{jn}\int_{\mathbb{R}^{n}}\psi_{0}(2^{j}(x-y))(\Delta_{j^{\prime}%
}f_{i}(y)-\dot{S}_{j-2}\Delta_{j^{\prime}}f_{i}(y))\Delta_{j^{\prime\prime}%
}\partial_{i}g_{k}(y)\,dy\\
&  =\sum_{i=1}^{n}\sum_{|j^{\prime}-j^{\prime\prime}|\leq1}\sum_{j^{\prime
}\geq j-3}2^{j}2^{jn}\int_{\mathbb{R}^{n}}\partial_{i}\psi_{0}(2^{j}%
(x-y))(\Delta_{j^{\prime}}f_{i}(y)-\dot{S}_{j-2}\Delta_{j^{\prime}}%
f_{i}(y))\Delta_{j^{\prime\prime}}g_{k}(y)\,dy\\
&  =\sum_{i=1}^{n}\sum_{|j^{\prime}-j^{\prime\prime}|\leq1}\sum_{j^{\prime
}\geq j-3}2^{j}\{(2^{jn}\partial_{i}\psi_{0}(2^{j}\cdot))\ast(\Delta
_{j^{\prime}}f_{i}-\dot{S}_{j-2}\Delta_{j^{\prime}}f_{i})\Delta_{j^{\prime
\prime}}g_{k})\}
\end{split}
\]
Then
\begin{equation}%
\begin{split}
\Vert IV_{1}\Vert_{L^{p}}  &  \leq C\sum_{i=1}^{n}\sum_{|j^{\prime}%
-j^{\prime\prime}|\leq1}\sum_{j^{\prime}\geq j-3}2^{j}(\Vert\Delta_{j^{\prime
}}f_{i}\Vert_{L^{p}}+\Vert\dot{S}_{j-2}\Delta_{j^{\prime}}f_{i}\Vert_{L^{p}%
})\Vert\Delta_{j^{\prime\prime}}g_{k}\Vert_{L^{\infty}}\\
&  \leq C\Vert g\Vert_{L^{p_{2}}}\sum_{j^{\prime}\geq j-3}2^{j}\Vert
\Delta_{j^{\prime}}f\Vert_{L^{p_{1}}}.
\end{split}
\label{estimate_IV_1}%
\end{equation}
On the other hand, note that
\[%
\begin{split}
IV_{2}  &  =\sum_{i=1}^{n}\sum_{|j^{\prime}-j^{\prime\prime}|\leq1}[\Delta
_{j}((\Delta_{j^{\prime}}\dot{S}_{j-2}f_{i})\Delta_{j^{\prime\prime}}%
\partial_{i}g_{k})-(\Delta_{j^{\prime}}\dot{S}_{j-2}f_{i})(\Delta
_{j^{\prime\prime}}\Delta_{j}\partial_{i}g_{k})]\\
&  =\sum_{i=1}^{n}\sum_{j-1\geq j^{\prime}\geq j-3}\sum_{|j^{\prime}%
-j^{\prime\prime}|\leq1}[\Delta_{j},\Delta_{j^{\prime}}\dot{S}_{j-2}%
f_{i}]\Delta_{j^{\prime\prime}}\partial_{i}g_{k}.
\end{split}
\]
Also, we can write
\[%
\begin{split}
\lbrack\Delta_{j},  &  \Delta_{j^{\prime}}\dot{S}_{j-2}f_{i}]\Delta
_{j^{\prime\prime}}\partial_{i}g_{k}\\
&  =2^{jn}\int_{\mathbb{R}^{n}}\psi_{0}(2^{j}(x-y))(\Delta_{j^{\prime}}\dot
{S}_{j-2}f_{i}(y)-\Delta_{j^{\prime}}\dot{S}_{j-2}f_{i}(x))\Delta
_{j^{\prime\prime}}\partial_{i}g_{k}(y)\,dy\\
&  =2^{j(n+1)}\int_{\mathbb{R}^{n}}\partial_{i}\psi_{0}(2^{j}(x-y))(\Delta
_{j^{\prime}}\dot{S}_{j-2}f_{i}(y)-\Delta_{j^{\prime}}\dot{S}_{j-2}%
f_{i}(x))\Delta_{j^{\prime\prime}}g_{k}(y)\,dy\\
&  =2^{j(n+1)}\int_{\mathbb{R}^{n}}\partial_{i}\psi_{0}(2^{j}(x-y))\int
_{0}^{1}((y-x)\cdot\nabla)(\Delta_{j^{\prime}}\dot{S}_{j-2}f_{i}%
(x+\tau(y-x))\,d\tau\Delta_{j^{\prime\prime}}g_{k}(y)\,dy\\
&  =\int_{\mathbb{R}^{n}}\partial_{i}\psi_{0}(z)\int_{0}^{1}(z\cdot
\nabla)(\Delta_{j^{\prime}}\dot{S}_{j-2}f_{i}(x-\tau2^{-j}z)\,d\tau
\Delta_{j^{\prime\prime}}g_{k}(x-2^{-j}z)\,dz.
\end{split}
\]
Hence,
\[%
\begin{split}
\Vert\lbrack\Delta_{j},\Delta_{j^{\prime}}\dot{S}_{j-2}f_{i}]\Delta
_{j^{\prime\prime}}\partial_{i}g_{k}\Vert_{L^{p}}  &  \leq C\sum_{m=1}%
^{n}\Vert\Delta_{j^{\prime}}\dot{S}_{j-2}\partial_{m}f_{i}\Vert_{L^{\infty}%
}\Vert\Delta_{j^{\prime\prime}}g_{k}\Vert_{L^{p}}\\
&  \leq C\Vert\nabla f\Vert_{L^{\infty}}\Vert\Delta_{j^{\prime\prime}}%
g_{k}\Vert_{L^{p}},
\end{split}
\]
and
\begin{equation}
\Vert IV_{2}\Vert_{L^{p}}\leq C\Vert\nabla f\Vert_{L^{r_{2}}}\sum
_{|j-j^{\prime}|\leq5}\Vert\Delta_{j^{\prime}}g_{k}\Vert_{L^{r_{1}}}.
\label{estimate_IV_2}%
\end{equation}
Summing up the estimates (\ref{estimate_I}), (\ref{estimate_II}),
(\ref{estimate_III}), (\ref{estimate_IV_1}) and (\ref{estimate_IV_2}), we
obtain
\[%
\begin{split}
\Vert &  (\dot{S}_{j-2}f\cdot\nabla)\Delta_{j}g-\Delta_{j}(f\cdot\nabla
)g\Vert_{L^{p}}\\
&  \leq C\Vert\nabla f\Vert_{L^{r_{2}}}\sum_{|j-j^{\prime}|\leq5}\Vert
\Delta_{j^{\prime}}g_{k}\Vert_{L^{r_{1}}}+C\Vert g\Vert_{L^{p_{2}}}\left(
\sum_{|j-j^{\prime}|\leq3}2^{j}\Vert\Delta_{j^{\prime}}f\Vert_{L^{p_{1}}}%
+\sum_{j^{\prime}\geq j-3}2^{j}\Vert\Delta_{j^{\prime}}f\Vert_{L^{p_{1}}%
}\right)  .
\end{split}
\]
Multiplying by $2^{js}$ and computing the $l^{q}(\mathbb{Z})$-norm, we can
estimate
\[%
\begin{split}
&  \left(  \sum_{j\in\mathbb{Z}}2^{jqs}\Vert(\dot{S}_{j-2}f\cdot\nabla
)\Delta_{j}g-\Delta_{j}(f\cdot\nabla)g\Vert_{L^{p}}^{q}\right)  ^{1/q}\\
\leq &  C\Vert\nabla f\Vert_{L^{r_{2}}}\left(  \sum_{j\in\mathbb{Z}}%
\sum_{|j-j^{\prime}|\leq5}2^{jqs}\Vert\Delta_{j^{\prime}}g_{k}\Vert_{L^{r_{1}%
}}^{q}\right)  ^{1/q}+C\Vert g\Vert_{L^{p_{2}}}\left(  \sum_{j\in\mathbb{Z}%
}\sum_{|j-j^{\prime}|\leq3}2^{j(s+1)q}\Vert\Delta_{j^{\prime}}f\Vert
_{L^{p_{1}}}^{q}\right)  ^{1/q}\\
&  +C\Vert g\Vert_{L^{p_{2}}}\left(  \sum_{j\in\mathbb{Z}}\sum_{j^{\prime}\geq
j-3}2^{j(s+1)q}\Vert\Delta_{j^{\prime}}f\Vert_{L^{p_{1}}}^{q}\right)  ^{1/q}\\
:=  &  K_{1}+K_{2}+K_{3}.
\end{split}
\]
Now, observing that%
\[
\sum_{j\in\mathbb{Z}}\sum_{|j-j^{\prime}|\leq5}2^{jqs}\Vert\Delta_{j^{\prime}%
}h\Vert_{L^{p}}^{q}=\sum_{k=-5}^{5}2^{-kqs}\sum_{j\in\mathbb{Z}}%
2^{(j+k)qs}\Vert\Delta_{j+k}h\Vert_{L^{p}}^{q}\leq C\sum_{j\in\mathbb{Z}}%
\Vert\Delta_{j}h\Vert_{L^{p}}^{q},
\]
we have that
\[
K_{1}\leq C\Vert\nabla f\Vert_{L^{r_{2}}}\Vert g\Vert_{\dot{B}_{r_{1},q}^{s}%
},
\]
and similarly $K_{2}\leq C\Vert g\Vert_{L^{p_{2}}}\Vert f\Vert_{\dot{B}%
_{p_{1},q}^{s+1}}.$

For $K_{3}$, note that
\[%
\begin{split}
K_{3}  &  =C\Vert g\Vert_{L^{p_{2}}}\left(  \sum_{j\in\mathbb{Z}}%
\sum_{j^{\prime}\geq j-3}2^{(j-j^{\prime})(s+1)q}(2^{j^{\prime}(s+1)}%
\Vert\Delta_{j^{\prime}}f\Vert_{L^{p_{1}}})^{q}\right)  ^{1/q}\\
&  =C\Vert g\Vert_{L^{p_{2}}}\left(  \sum_{k\geq-3}2^{-kq(s+1)}\sum
_{j\in\mathbb{Z}}(2^{(j+k)(s+1)}\Vert\Delta_{j+k}f\Vert_{L^{p_{1}}}%
)^{q}\right)  ^{1/q}\\
&  =C\Vert g\Vert_{L^{p_{2}}}\left(  \sum_{k\geq-3}2^{-kq(s+1)}\right)
^{1/q}\left(  \sum_{j\in\mathbb{Z}}2^{jq(s+1)}\Vert\Delta_{j}f\Vert_{L^{p_{1}%
}}^{q}\right)  ^{1/q}\\
&  \leq C\Vert g\Vert_{L^{p_{2}}}\Vert f\Vert_{\dot{B}_{p_{1},q}^{s+1}}.
\end{split}
\]
Therefore,
\[
\left(  \sum_{j\in\mathbb{Z}}2^{sjq}\left\Vert (\dot{S}_{j-2}f\cdot
\nabla)\Delta_{j}g-\Delta_{j}(f\cdot\nabla)g\right\Vert _{L^{p}}^{q}\right)
^{1/q}\leq C\left(  \Vert\nabla f\Vert_{L^{r_{2}}}\Vert g\Vert_{\dot{B}%
_{r_{1},q}^{s}}+\Vert g\Vert_{L^{p_{2}}}\Vert f\Vert_{\dot{B}_{p_{1},q}^{s+1}%
}\right)  .
\]
This completes the proof of $(ii)$. \begin{flushright}$\blacksquare$%
\end{flushright}

\section{An approximate linear iteration problem and local-in-time
solvability}

\label{Localregularizationsolution} In order to prove the local existence to
(\ref{Bouss1}), we consider the approximate linear iteration problem
\begin{equation}
\left\{
\begin{split}
&  \partial_{t}\omega_{n+1}+(u_{n}\cdot\nabla)\omega_{n+1}=\kappa\partial
_{1}\rho_{n+1}\ \ \text{ in }\ \ \mathbb{R}^{2}\times(0,\infty),\\
&  \partial_{t}\rho_{n+1}+(u_{n}\cdot\nabla)\rho_{n+1}=\kappa u_{2,n+1}%
\ \ \text{ in }\ \ \mathbb{R}^{2}\times(0,\infty),\\
&  u_{n+1}=\nabla^{\perp}(-\Delta)^{-1}\omega_{n+1}\ \ \text{ in
}\ \ \mathbb{R}^{2}\times(0,\infty),\\
&  \omega_{n+1}|_{t=0}=S_{n+2}\omega_{0},\ \rho_{n+1}|_{t=0}=S_{n+2}\rho
_{0}\ \ \text{ in }\mathbb{R}^{2}.
\end{split}
\right.  \label{BCcurl2approx}%
\end{equation}
From (\ref{BCcurl2approx}), we provide uniform estimates for the sequence
$\left\{  (\omega_{n},\rho_{n})\right\}  _{n\in\mathbb{N}}$ and then obtain a
solution for (\ref{Bouss1}).

\textbf{Uniform estimates.} Applying $\Delta_{j}$ in (\ref{BCcurl2approx}),
taking the product with $\Delta_{j}\omega_{n+1}$ in $\dot{H}^{-1}$ and the
product with $\Delta_{j}\rho_{n+1}$ in $L^{2}$ in the first and second
equations, respectively, and using the divergence-free condition $\nabla
\cdot\Delta_{j}u_{n}=0$, we obtain that
\[%
\begin{split}
\langle\Delta_{j}\partial_{t}\omega_{n+1},\Delta_{j}\omega_{n+1}\rangle
_{\dot{H}^{-1}}=  &  \langle\lbrack u_{n}\cdot\nabla,\Delta_{j}]\omega
_{n+1},\Delta_{j}\omega_{n+1}\rangle_{\dot{H}^{-1}}+\kappa\langle\Delta
_{j}\partial_{1}\rho_{n+1},\Delta_{j}\omega_{n+1}\rangle_{\dot{H}^{-1}},\\
\langle\Delta_{j}\partial_{t}\rho_{n+1},\Delta_{j}\rho_{n+1}\rangle_{L^{2}}=
&  \langle\lbrack u_{n}\cdot\nabla,\Delta_{j}]\rho_{n+1},\Delta_{j}\rho
_{n+1}\rangle_{L^{2}}+\kappa\langle\Delta_{j}u_{2,n+1},\Delta_{j}\rho
_{n+1}\rangle_{L^{2}}.
\end{split}
\]
Adding the two above inequalities and using the properties
\begin{equation}%
\begin{split}
&  \langle(u_{n}\cdot\nabla)\Lambda^{-1}\Delta_{j}\omega_{n+1},\Lambda
^{-1}\Delta_{j}\omega_{n+1}\rangle_{L^{2}}=\langle(u_{n}\cdot\nabla)\Delta
_{j}\rho_{n+1},\Delta_{j}\rho_{n+1}\rangle_{L^{2}}=0,\\
&  \langle\Delta_{j}\partial_{1}\rho_{n+1},\Delta_{j}\omega_{n+1}\rangle
_{\dot{H}^{-1}}+\langle\Delta_{j}u_{2,n+1},\Delta_{j}\rho_{n+1}\rangle_{L^{2}%
}=0,
\end{split}
\label{cancelation}%
\end{equation}
we arrive at
\[%
\begin{split}
\frac{1}{2}\frac{d}{dt}\left(  \Vert\Delta_{j}\omega_{n+1}\Vert_{\dot{H}^{-1}%
}^{2}+\Vert\Delta_{j}\rho_{n+1}\Vert_{L^{2}}^{2}\right)  \leq\,  &  \left\Vert
[u_{n}\cdot\nabla,\Lambda^{-1}\Delta_{j}]\omega_{n+1}\right\Vert _{L^{2}%
}\left\Vert \Delta_{j}\omega_{n+1}\right\Vert _{\dot{H}^{-1}}\\
&  +\left\Vert [u_{n}\cdot\nabla,\Delta_{j}]\rho_{n+1}\right\Vert _{L^{2}%
}\left\Vert \Delta_{j}\rho_{n+1}\right\Vert _{L^{2}}.
\end{split}
\]
Integrating over $(0,t)$, it follows that
\[%
\begin{split}
\frac{1}{2}  &  \left(  \Vert\Delta_{j}\omega_{n+1}(t)\Vert_{\dot{H}^{-1}}%
^{2}+\Vert\Delta_{j}\rho_{n+1}(t)\Vert_{L^{2}}^{2}\right)  \leq\frac{1}%
{2}\left(  \Vert\Delta_{j}\omega_{n+1}(0)\Vert_{\dot{H}^{-1}}^{2}+\Vert
\Delta_{j}\rho_{n+1}(0)\Vert_{L^{2}}^{2}\right) \\
&  +\int_{0}^{t}(\left\Vert [u_{n}(\tau)\cdot\nabla,\Lambda^{-1}\Delta
_{j}]\omega_{n+1}(\tau)\right\Vert _{L^{2}}+\left\Vert [u_{n}(\tau)\cdot
\nabla,\Delta_{j}]\rho_{n+1}(\tau)\right\Vert _{L^{2}})\left(  \left\Vert
\Delta_{j}\omega_{n+1}(\tau)\right\Vert _{\dot{H}^{-1}}^{2}+\left\Vert
\Delta_{j}\rho_{n+1}(\tau)\right\Vert _{L^{2}}^{2}\right)  ^{\frac{1}{2}%
}\,d\tau.
\end{split}
\]
Thus, by Gr\"{o}nwall's inequality (see Proposition 1.2 in \cite[page
24]{Barbu}), we have
\[%
\begin{split}
\Vert\Delta_{j}\omega_{n+1}(t)\Vert_{\dot{H}^{-1}}+\Vert\Delta_{j}\rho
_{n+1}(t)\Vert_{L^{2}}\leq &  \Vert\Delta_{j}\omega_{n+1}(0)\Vert_{\dot
{H}^{-1}}+\Vert\Delta_{j}\rho_{n+1}(0)\Vert_{L^{2}}\\
&  +\int_{0}^{t}\left\Vert [u_{n}(\tau)\cdot\nabla,\Lambda^{-1}\Delta
_{j}]\omega_{n+1}(\tau)\right\Vert _{L^{2}}+\left\Vert [u_{n}(\tau)\cdot
\nabla,\Delta_{j}]\rho_{n+1}(\tau)\right\Vert _{L^{2}}\,d\tau.
\end{split}
\]
Multiplying by $2^{sj}$, applying the $l^{q}(\mathbb{Z})$-norm and Lemma
\ref{inequalitybernstein}, it follows that
\[%
\begin{split}
\Vert &  \omega_{n+1}(t)\Vert_{\dot{B}_{2,q}^{s-1}}+\Vert\rho_{n+1}%
(t)\Vert_{\dot{B}_{2,q}^{s}}\leq C\Vert\omega_{n+1}(0)\Vert_{\dot{B}%
_{2,q}^{s-1}}+C\Vert\rho_{n+1}(0)\Vert_{\dot{B}_{2,q}^{s}}\\
&  +C\int_{0}^{t}\left(  \sum_{j\in\mathbb{Z}}2^{sjq}\left\Vert [u_{n}%
(\tau)\cdot\nabla,\Lambda^{-1}\Delta_{j}]\omega_{n+1}(\tau)\right\Vert
_{L^{2}}^{q}\right)  ^{\frac{1}{q}}+\left(  \sum_{j\in\mathbb{Z}}%
2^{sjq}\left\Vert [u_{n}(\tau)\cdot\nabla,\Delta_{j}]\rho_{n+1}(\tau
)\right\Vert _{L^{2}}^{q}\right)  ^{\frac{1}{q}}\,d\tau.
\end{split}
\]
From Remark \ref{remark2} and Lemma \ref{2.4}, we have
\[%
\begin{split}
\left(  \sum_{j\in\mathbb{Z}}2^{sjq}\left\Vert [u_{n}\cdot\nabla,\Lambda
^{-1}\Delta_{j}]\omega_{n+1}\right\Vert _{L^{2}}^{q}\right)  ^{\frac{1}{q}}
&  \leq C(\left\Vert \nabla u_{n}\right\Vert _{L^{\infty}}\left\Vert
\omega_{n+1}\right\Vert _{\dot{B}_{2,q}^{s-1}}+\left\Vert \omega
_{n+1}\right\Vert _{L^{2}}\left\Vert u_{n}\right\Vert _{\dot{B}_{\infty,q}%
^{s}})\\
&  \leq C\left\Vert \omega_{n}\right\Vert _{\dot{B}_{2,q}^{s-1}\cap\dot
{H}^{-1}}\left\Vert \omega_{n+1}\right\Vert _{\dot{B}_{2,q}^{s-1}\cap\dot
{H}^{-1}},\\
\left(  \sum_{j\in\mathbb{Z}}2^{sjq}\left\Vert [u_{n}\cdot\nabla,\Delta
_{j}]\rho_{n+1}\right\Vert _{L^{2}}^{q}\right)  ^{\frac{1}{q}}  &  \leq
C(\left\Vert \nabla u_{n}\right\Vert _{L^{\infty}}\left\Vert \rho
_{n+1}\right\Vert _{\dot{B}_{2,q}^{s}}+\left\Vert \nabla\rho_{n+1}\right\Vert
_{L^{\infty}}\left\Vert u_{n}\right\Vert _{\dot{B}_{2,q}^{s}})\\
&  \leq C\left\Vert \omega_{n}\right\Vert _{\dot{B}_{2,q}^{s-1}\cap\dot
{H}^{-1}}\left\Vert \rho_{n+1}\right\Vert _{B_{2,q}^{s}}.
\end{split}
\]
Then,
\begin{equation}%
\begin{split}
\Vert\omega_{n+1}(t)\Vert_{\dot{B}_{2,q}^{s-1}}+\Vert\rho_{n+1}(t)\Vert
_{\dot{B}_{2,q}^{s}}\leq &  \,C\Vert\omega_{n+1}(0)\Vert_{\dot{B}_{2,q}^{s-1}%
}+C\Vert\rho_{n+1}(0)\Vert_{\dot{B}_{2,q}^{s}}\\
&  +C\int_{0}^{t}\left\Vert \omega_{n}\right\Vert _{\dot{B}_{2,q}^{s-1}%
\cap\dot{H}^{-1}}(\left\Vert \omega_{n+1}\right\Vert _{\dot{B}_{2,q}^{s-1}%
\cap\dot{H}^{-1}}+\left\Vert \rho_{n+1}\right\Vert _{B_{2,q}^{s}})\,d\tau.
\end{split}
\label{Uniform_hom_besov}%
\end{equation}
On the other hand, taking the product with $\omega_{n+1}$ in $\dot{H}^{-1}$
and the product with $\rho_{n+1}$ in $L^{2}$ in the first and second equations
of (\ref{BCcurl2approx}), respectively, and using the divergence-free
condition $\nabla\cdot\Delta_{j}u_{n}=0$, we obtain that
\[%
\begin{split}
\frac{1}{2}\frac{d}{dt}\left(  \left\Vert \omega_{n+1}\right\Vert _{\dot
{H}^{-1}}^{2}+\left\Vert \rho_{n+1}\right\Vert _{L^{2}}^{2}\right)   &
\leq\left\Vert \Lambda^{-1}\nabla\cdot(u_{n}\otimes\omega_{n+1})\right\Vert
_{L^{2}}\left\Vert \omega_{n+1}\right\Vert _{\dot{H}^{-1}}\\
&  \leq C\left\Vert u_{n}\right\Vert _{L^{2}}\left\Vert \omega_{n+1}%
\right\Vert _{L^{\infty}}\left\Vert \omega_{n+1}\right\Vert _{\dot{H}^{-1}}\\
&  \leq C\left\Vert \omega_{n}\right\Vert _{\dot{B}_{2,q}^{s-1}\cap\dot
{H}^{-1}}\left\Vert \omega_{n+1}\right\Vert _{\dot{B}_{2,q}^{s-1}\cap\dot
{H}^{-1}}\left\Vert \omega_{n+1}\right\Vert _{\dot{H}^{-1}}.
\end{split}
\]
Here, we have used the equality $u_{n}=\nabla^{\perp}(-\Delta)^{-1}\omega_{n}%
$, Lemma \ref{inequalitybernstein}, Remark \ref{remark1}, the embedding
$\dot{B}_{2,q}^{s-1}\cap\dot{H}^{-1}\hookrightarrow B_{2,q}^{s-1}$ and the
property (\ref{cancelation}). Now, we integrate over $(0,t)$ and we use a
Gr\"{o}nwall-type inequality to get
\begin{equation}%
\begin{split}
\left\Vert \omega_{n+1}(t)\right\Vert _{\dot{H}^{-1}}+\left\Vert \rho
_{n+1}(t)\right\Vert _{L^{2}}\leq\,  &  C\left\Vert \omega_{n+1}(0)\right\Vert
_{\dot{H}^{-1}}+C\left\Vert \rho_{n+1}(0)\right\Vert _{L^{2}}\\
&  +C\int_{0}^{t}\left\Vert \omega_{n}(\tau)\right\Vert _{\dot{B}_{2,q}%
^{s-1}\cap\dot{H}^{-1}}\left\Vert \omega_{n+1}(\tau)\right\Vert _{\dot
{B}_{2,q}^{s-1}\cap\dot{H}^{-1}}\,d\tau.
\end{split}
\label{Umiform_HL2}%
\end{equation}
Combining (\ref{Uniform_hom_besov}) and (\ref{Umiform_HL2}), we obtain
\[%
\begin{split}
\left\Vert \omega_{n+1}(t)\right\Vert _{\dot{B}_{2,q}^{s-1}\cap\dot{H}^{-1}}
&  +\left\Vert \rho_{n+1}(t)\right\Vert _{B_{2,q}^{s}}\leq\left\Vert
\omega_{n+1}(0)\right\Vert _{\dot{B}_{2,q}^{s-1}\cap\dot{H}^{-1}}+\left\Vert
\rho_{n+1}(0)\right\Vert _{B_{2,q}^{s}}\\
&  +C\int_{0}^{t}\left\Vert \omega_{n}(\tau)\right\Vert _{\dot{B}_{2,q}%
^{s-1}\cap\dot{H}^{-1}}(\left\Vert \omega_{n+1}(\tau)\right\Vert _{\dot
{B}_{2,q}^{s-1}\cap\dot{H}^{-1}}+\left\Vert \rho_{n+1}(\tau)\right\Vert
_{B_{2,q}^{s}})\,d\tau.
\end{split}
\]
Doing $A_{n+1}(t):=\left\Vert \omega_{n+1}(t)\right\Vert _{\dot{B}_{2,q}%
^{s-1}\cap\dot{H}^{-1}}+\left\Vert \rho_{n+1}(t)\right\Vert _{B_{2,q}^{s}}$,
we have
\[
A_{n+1}(t)\leq A_{n+1}(0)+C\int_{0}^{t}A_{n}(\tau)A_{n+1}(\tau)\ d\tau.
\]
Since $\left\Vert \omega_{n+1}(0)\right\Vert _{\dot{B}_{2,q}^{s-1}\cap\dot
{H}^{-1}}\leq C\left\Vert \omega_{0}\right\Vert _{\dot{B}_{2,q}^{s-1}\cap
\dot{H}^{-1}}$ and $\left\Vert \rho_{n+1}(0)\right\Vert _{B_{2,q}^{s}}\leq
C\left\Vert \rho_{0}\right\Vert _{B_{2,q}^{s}}$, Gr\"{o}nwall's inequality
yields
\[
A_{n+1}(t)\leq C_{0}A_{0}\exp\left(  C_{1}\int_{0}^{t}A_{n}(\tau
)\ d\tau\right)  ,
\]
where $A_{0}=\left\Vert \omega_{0}\right\Vert _{\dot{B}_{2,q}^{s-1}\cap\dot
{H}^{-1}}+\left\Vert \rho_{0}\right\Vert _{B_{2,q}^{s}}$ and the constants
$C_{0},C_{1}>0$ are independent of $n.$ We wish to show that there exist $P>0$
and $T>0$ satisfying
\begin{equation}
A_{n+1}(t)\leq PA_{0},\ \ \text{ for all }t\in(0,T)\text{ and }n\in\mathbb{N}.
\label{unif_est}%
\end{equation}
In fact, for $n=0$, note that
\[
A_{1}(t)\leq C_{0}A_{0}\exp\left(  C_{1}tA_{0}\right)  \leq PA_{0},\text{ for
all }t\in(0,T_{1}),
\]
where $T_{1}:=\frac{1}{C_{1}A_{0}}\log\left(  \frac{P}{C_{0}}\right)  $.
Similarly, for $n=1$ we have
\[
A_{2}(t)\leq C_{0}A_{0}\exp\left(  C_{1}\int_{0}^{t}A_{1}(\tau)\,d\tau\right)
\leq C_{0}A_{0}\exp\left(  C_{1}tPA_{0}\right)  \leq PA_{0},
\]
for all $t\in(0,T_{2})$, where $T_{2}:=\frac{1}{C_{1}PA_{0}}\log\left(
\frac{P}{C_{0}}\right)  $. Denoting $T_{3}:=\min\{T_{1},T_{2}\}$ and making
the same calculations for $n=2$, we arrive at
\[
A_{3}(t)\leq PA_{0},\text{ for all }t\in(0,T_{3}).
\]
Thus, continuing in the same way, we obtain (\ref{unif_est}) by induction.

\textbf{Continuity of the sequence.} Our intent now is to show that the
sequences $\{\omega_{n}\}_{n\in\mathbb{N}}$ and $\{\rho_{n}\}_{n\in\mathbb{N}%
}$ belong to $C([0,T];\dot{B}_{2,q}^{s-1}(\mathbb{R}^{2})\cap\dot{H}%
^{-1}(\mathbb{R}^{2}))$ and $C([0,T];B_{2,q}^{s}(\mathbb{R}^{2}))$,
respectively. For that, considering the equality $(f\cdot\nabla)g=\nabla
\cdot(f\otimes g)$ for all divergence free vector fields $f$, Remark
\ref{remark1}, Lemma \ref{inequalityholder}, the embedding $\dot{B}%
_{2,q}^{s-1}\cap\dot{H}^{-1}\hookrightarrow B_{2,q}^{s-1}$, the following
estimates
\[%
\begin{split}
\Vert(u_{n}\cdot\nabla)\omega_{n+1}\Vert_{\dot{B}_{2,q}^{s-2}\cap\dot{H}%
^{-1}}  &  \leq C\Vert u_{n}\otimes\omega_{n+1}\Vert_{\dot{B}_{2,q}^{s-1}\cap
L^{2}}\leq C\Vert\omega_{n}\Vert_{\dot{B}_{2,q}^{s-1}\cap\dot{H}^{-1}}%
\Vert\omega_{n+1}\Vert_{\dot{B}_{2,q}^{s-1}\cap\dot{H}^{-1}},\\
\Vert(u_{n}\cdot\nabla)\rho_{n+1}\Vert_{B_{2,q}^{s-1}}  &  \leq C\Vert
\omega_{n}\Vert_{\dot{B}_{2,q}^{s-1}\cap\dot{H}^{-1}}\Vert\rho_{n+1}%
\Vert_{B_{2,q}^{s}},
\end{split}
\]
estimate (\ref{unif_est}), and the two first equations of (\ref{BCcurl2approx}%
), we have that $\partial_{t}\omega_{n+1}\in L^{\infty}(0,T;\dot{B}%
_{2,q}^{s-1}(\mathbb{R}^{2})\cap\dot{H}^{-1}(\mathbb{R}^{2}))$ and
$\partial_{t}\rho_{n+1}\in L^{\infty}(0,T;B_{2,q}^{s-1}(\mathbb{R}^{2}))$.
Thus,
\[%
\begin{split}
\omega_{n+1}  &  \in W^{1,\infty}([0,T];\dot{B}_{2,q}^{s-1}(\mathbb{R}%
^{2})\cap\dot{H}^{-1}(\mathbb{R}^{2}))\subset C([0,T];\dot{B}_{2,q}%
^{s-1}(\mathbb{R}^{2})\cap\dot{H}^{-1}(\mathbb{R}^{2})),\\
\rho_{n+1}  &  \in W^{1,\infty}([0,T];B_{2,q}^{s-1}(\mathbb{R}^{2}))\subset
C([0,T];B_{2,q}^{s-1}(\mathbb{R}^{2})).
\end{split}
\]
For $k\in\mathbb{N}$ and $n$ fixed, we denote $y_{k}^{n}:=\dot{S}_{k}%
\omega_{n+1}$ and $z_{k}^{n}:=S_{k}\rho_{n+1}$. We claim that $y_{k}%
^{n}\rightarrow\omega_{n+1}$ in $L^{\infty}(0,T;\dot{B}_{2,q}^{s-1}%
(\mathbb{R}^{2})\cap\dot{H}^{-1}(\mathbb{R}^{2}))$ and $z_{k}^{n}%
\rightarrow\rho_{n+1}$ in $L^{\infty}(0,T;B_{2,q}^{s}(\mathbb{R}^{2}))$,
respectively, as $k\rightarrow\infty$. Using the Littlewood-Paley operators,
we can also write
\[
\partial_{t}\Delta_{j}\omega_{n+1}+(\dot{S}_{j-2}u_{n}\cdot\nabla)\Delta
_{j}\omega_{n+1}=\kappa\Delta_{j}\partial_{1}\rho_{n+1}+(\dot{S}_{j-2}%
u_{n}\cdot\nabla)\Delta_{j}\omega_{n+1}-\Delta_{j}(u_{n}\cdot\nabla
)\omega_{n+1},
\]
for each $j\in\mathbb{N}$. Since $\Delta_{j}\omega_{n+1}$ and $\Delta_{j}%
\rho_{n+1}$ are absolutely continuous functions from $[0,T]$ to $L^{2}%
(\mathbb{R}^{2})$ and $\nabla\cdot S_{j-2}u_{n}=0$, we obtain that
\[%
\begin{split}
\Vert\Delta_{j}\omega_{n+1}(t)\Vert_{\dot{H}^{-1}}\leq\,  &  \Vert\Delta
_{j}\omega_{n+1}(0)\Vert_{\dot{H}^{-1}}+\left\vert \kappa\right\vert \int
_{0}^{t}\Vert\Delta_{j}\partial_{1}\rho_{n+1}(\tau)\Vert_{\dot{H}^{-1}}%
\,d\tau\\
&  +\int_{0}^{t}\Vert(\dot{S}_{j-2}u_{n}(\tau)\cdot\nabla)\Delta_{j}%
\omega_{n+1}(\tau)-\Delta_{j}(u_{n}(\tau)\cdot\nabla)\omega_{n+1}(\tau
)\Vert_{\dot{H}^{-1}}\,d\tau,\\
\Vert\Delta_{j}\rho_{n+1}(t)\Vert_{L^{2}}\leq\,  &  \Vert\Delta_{j}\rho
_{n+1}(0)\Vert_{L^{2}}+|\kappa|\int_{0}^{t}\Vert\Delta_{j}u_{2,n+1}(\tau
)\Vert_{L^{2}}\,d\tau+\int_{0}^{t}\Vert\Delta_{j}(u_{n}(\tau)\cdot\nabla
)\rho_{n+1}\Vert_{L^{2}}\,d\tau.
\end{split}
\]
It follows that
\[%
\begin{split}
\Vert\omega_{n+1}(t)  &  -y_{k}^{n}(t)\Vert_{\dot{B}_{2,q}^{s-1}}\leq C\left(
\sum_{j> k}2^{j(s-1)q}\Vert\Delta_{j}\omega_{n+1}(t)\Vert_{L^{2}}^{q}\right)
^{1/q}\\
\leq\,  &  C\left(  \sum_{j> k}2^{j(s-1)q}\Vert\Delta_{j}\omega_{n+1}%
(0)\Vert_{L^{2}}^{q}\right)  ^{1/q}+C|\kappa|\int_{0}^{t}\left(  \sum_{j>
k}2^{j(s-1)q}\Vert\Delta_{j}\partial_{1}\rho_{n+1}(\tau)\Vert_{L^{2}}%
^{q}\right)  ^{1/q}\,d\tau\\
&  +C\int_{0}^{t}\Bigl(\sum_{j> k}2^{j(s-1)q}\Vert(\dot{S}_{j-2}u_{n}%
(\tau)\cdot\nabla)\Delta_{j}\omega_{n+1}(\tau)-\Delta_{j}(u_{n}(\tau
)\cdot\nabla)\omega_{n+1}(\tau)\Vert_{L^{2}}^{q}\Bigr) ^{1/q}\,d\tau,
\end{split}
\]
\[%
\begin{split}
\Vert\omega_{n+1}(t)-y_{k}^{n}(t)\Vert_{\dot{H}^{-1}}\leq C\Vert\Delta
_{j}\omega_{n+1}(t)\Vert_{\dot{H}^{-1}}\leq\,  &  C\Vert\Delta_{j}\omega
_{n+1}(0)\Vert_{\dot{H}^{-1}}+C|\kappa|\int_{0}^{t}\Vert\Delta_{j}\partial
_{1}\rho_{n+1}(\tau)\Vert_{\dot{H}^{-1}}\,d\tau\\
&  +C\int_{0}^{t}\Vert\Delta_{j}(u_{n}(\tau)\cdot\nabla)\omega_{n+1}%
(\tau)\Vert_{\dot{H}^{-1}}\,d\tau,
\end{split}
\]
\[%
\begin{split}
\Vert\rho_{n+1}(t)-z_{k}^{n}(t)\Vert_{B_{2,q}^{s}}\leq\,  &  C\left(  \sum_{j>
k}2^{jsq}\Vert\Delta_{j}\rho_{n+1}(t)\Vert_{L^{2}}^{q}\right)  ^{1/q}\\
\leq\,  &  C\left(  \sum_{j> k}2^{jsq}\Vert\Delta_{j}\rho_{n+1}(0)\Vert
_{L^{2}}^{q}\right)  ^{1/q}+C|\kappa|\int_{0}^{t}\left(  \sum_{j> k}%
2^{jsq}\Vert\Delta_{j}u_{2,n+1}(\tau)\Vert_{L^{2}}^{q}\right)  ^{1/q}\,d\tau\\
&  +C\int_{0}^{t}\left(  \sum_{j> k}2^{jsq}\Vert\Delta_{j}(u_{n}(\tau
)\cdot\nabla)\rho_{n+1}(\tau)\Vert_{L^{2}}^{q}\right)  ^{1/q}\,d\tau.
\end{split}
\]
As $\omega_{n+1}(t)\in\dot{B}_{2,q}^{s-1}(\mathbb{R}^{2})\cap\dot{H}%
^{-1}(\mathbb{R}^{2})$ and $\rho_{n+1}(t)\in B_{2,q}^{s}(\mathbb{R}^{2})$, by
Lemmas \ref{inequalityholder} and \ref{estimate_type_commutator}, the R.H.S.
of the three above estimates go to zero as $k\rightarrow\infty$, and then the
desired claim follows. Moreover, we get
\begin{equation}%
\begin{split}
\Vert y_{k}^{n}(t^{\prime})  &  -y_{k}^{n}(t)\Vert_{\dot{B}_{2,q}^{s-1}%
\cap\dot{H}^{-1}} =\Vert\dot{S}_{k}(\omega_{n+1}(t^{\prime})-\omega
_{n+1}(t))\Vert_{\dot{B}_{2,q}^{s-1}\cap\dot{H}^{-1}}\\
&  \leq\left(  \sum_{j\leq k+1}2^{(s-1)jq}\Vert\Delta_{j}(\omega
_{n+1}(t^{\prime})-\omega_{n+1}(t))\Vert_{L^{2}}^{q}\right)  ^{1/q} +
\sum_{j\leq k+1}\Vert\Delta_{j}(\omega_{n+1}(t^{\prime})-\omega_{n+1}%
(t))\Vert_{\dot{H}^{-1}}\\
&  \leq C2^{k+1}\Vert\omega_{n+1}(t^{\prime})-\omega_{n+1}(t)\Vert_{\dot
{B}_{2,q}^{s-2}\cap\dot{H}^{-1}},\\
\Vert z_{k}^{n}(t^{\prime})  &  -z_{k}^{n}(t)\Vert_{B_{2,q}^{s}} =\Vert
S_{k}(\rho_{n+1}(t^{\prime})-\rho_{n+1}(t))\Vert_{B_{2,q}^{s}}\leq\left(
\sum_{j=-1}^{k+1}2^{sjq}\Vert\Delta_{j}(\rho_{n+1}(t^{\prime})-\rho
_{n+1}(t))\Vert_{L^{2}}^{q}\right)  ^{1/q}\\
&  \leq C2^{k+1}\Vert\rho_{n+1}(t^{\prime})-\rho_{n+1}(t)\Vert_{B_{2,q}^{s-1}%
}.
\end{split}
\label{aux-cont-7000}%
\end{equation}
Thus, $\{y_{k}^{n}\}_{k\in\mathbb{N}}\subset C([0,T];\dot{B}_{2,q}%
^{s-1}(\mathbb{R}^{2})\cap\dot{H}^{-1}(\mathbb{R}^{2}))$ and $\{z_{k}%
^{n}\}_{k\in\mathbb{N}}\subset C([0,T];B_{2,q}^{s}(\mathbb{R}^{2}))$. Then
\begin{equation}
\{\omega_{n}\}_{n\in\mathbb{N}}\subset C([0,T];\dot{B}_{2,q}^{s-1}%
(\mathbb{R}^{2})\cap\dot{H}^{-1}(\mathbb{R}^{2}))\ \text{ and }\ \{\rho
_{n}\}_{n\in\mathbb{N}}\subset C([0,T];B_{2,q}^{s}(\mathbb{R}^{2})).
\label{continuity_seq}%
\end{equation}

\textbf{Convergence and local solution.} Now, we show that $\{\omega
_{n}\}_{n\in\mathbb{N}}$ and $\{\rho_{n}\}_{n\in\mathbb{N}}$ converge in
$C([0,T];\dot{B}_{2,q}^{s-2}(\mathbb{R}^{2})\cap\dot{H}^{-1}(\mathbb{R}^{2}))$
and $C([0,T];B_{2,q}^{s-1}(\mathbb{R}^{2}))$, for some $T>0$, respectively.
For that, we consider the following system%
\begin{equation}
\left\{
\begin{split}
&  \partial_{t}\overline{\omega_{n+1}}+(\overline{u_{n}}\cdot\nabla
)\omega_{n+1}+(u_{n-1}\cdot\nabla)\overline{\omega_{n+1}}=\kappa\partial
_{1}\overline{\rho_{n+1}}\ \ \text{ in }\ \ \mathbb{R}^{2}\times(0,\infty),\\
&  \partial_{t}\overline{\rho_{n+1}}+(\overline{u_{n}}\cdot\nabla)\rho
_{n+1}+(u_{n-1}\cdot\nabla)\overline{\rho_{n+1}}=\kappa\overline{u_{2,n+1}%
}\ \ \text{ in }\ \ \mathbb{R}^{2}\times(0,\infty),\\
&  \overline{u_{n+1}}=\nabla^{\perp}(-\Delta)^{-1}\overline{\omega_{n+1}%
}\ \ \text{ in }\ \ \mathbb{R}^{2}\times(0,\infty),\\
&  \overline{\omega_{n+1}}\mid_{t=0}=\Delta_{n+1}\omega_{0},\ \overline
{\rho_{n+1}}\mid_{t=0}=\Delta_{n+1}\rho_{0}\ \ \text{ in }\mathbb{R}^{2},
\end{split}
\right.  \label{BCcurl2approx_dif}%
\end{equation}
where $\overline{\omega_{n+1}}:=\omega_{n+1}-\omega_{n}$ and $\overline
{\rho_{n+1}}:=\rho_{n+1}-\rho_{n}$.

We take the $\dot{H}^{-1}$-product with $\overline{\omega_{n+1}}$ and the
$L^{2}$-product with $\overline{\rho_{n+1}}$ in (\ref{BCcurl2approx_dif}) to
obtain
\[%
\begin{split}
\langle\partial_{t}\overline{\omega_{n+1}},\overline{\omega_{n+1}}%
\rangle_{\dot{H}^{-1}}+\langle(\overline{u_{n}}\cdot\nabla)\omega
_{n+1},\overline{\omega_{n+1}}\rangle_{\dot{H}^{-1}}+\langle(u_{n-1}%
\cdot\nabla)\overline{\omega_{n+1}},\overline{\omega_{n+1}}\rangle_{\dot
{H}^{-1}}  &  =\langle\kappa\partial_{1}\overline{\rho_{n+1}},\overline
{\omega_{n+1}}\rangle_{\dot{H}^{-1}},\\
\langle\partial_{t}\overline{\rho_{n+1}},\overline{\rho_{n+1}}\rangle_{L^{2}%
}+\langle(\overline{u_{n}}\cdot\nabla)\rho_{n+1},\overline{\rho_{n+1}}%
\rangle_{L^{2}}  &  =\langle\kappa\overline{u_{2,n+1}},\overline{\rho_{n+1}%
}\rangle_{L^{2}}.
\end{split}
\]
Next, adding the above two equalities, using the property $\langle\partial
_{1}\overline{\rho_{n+1}},\overline{\omega_{n+1}}\rangle_{\dot{H}^{-1}%
}+\langle\overline{u_{2,n+1}},\overline{\rho_{n+1}}\rangle_{L^{2}}=0$, the
equality $(f\cdot\nabla)g=\nabla\cdot(f\otimes g)$ for all divergence free
vector fields $f$, Cauchy-Schwarz and H\"{o}lder inequalities, Lemma
\ref{inequalitybernstein}, Remark \ref{remark1}, the embedding $\dot{B}%
_{2,q}^{s-2}\cap\dot{H}^{-1}\hookrightarrow L^{2}$ and the estimate
(\ref{unif_est}) to arrive at
\[%
\begin{split}
\frac{1}{2}\frac{d}{dt}  &  \left(  \Vert\overline{\omega_{n+1}}\Vert_{\dot
{H}^{-1}}^{2}+\Vert\overline{\rho_{n+1}}\Vert_{L^{2}}^{2}\right) \\
\leq\,  &  \Vert(\overline{u_{n}}\cdot\nabla)\omega_{n+1}\Vert_{\dot{H}^{-1}%
}\Vert\overline{\omega_{n+1}}\Vert_{\dot{H}^{-1}}+\Vert(u_{n-1}\cdot
\nabla)\overline{\omega_{n+1}}\Vert_{\dot{H}^{-1}}\Vert\overline{\omega_{n+1}%
}\Vert_{\dot{H}^{-1}}+\Vert(\overline{u_{n}}\cdot\nabla)\rho_{n+1}\Vert
_{L^{2}}\Vert\overline{\rho_{n+1}}\Vert_{L^{2}}\\
\leq\,  &  C\Vert\overline{u_{n}}\Vert_{L^{2}}\Vert\omega_{n+1}\Vert
_{L^{\infty}}\Vert\overline{\omega_{n+1}}\Vert_{\dot{H}^{-1}}+C\Vert
u_{n-1}\Vert_{L^{\infty}}\Vert\overline{\omega_{n+1}}\Vert_{L^{2}}%
\Vert\overline{\omega_{n+1}}\Vert_{\dot{H}^{-1}}+C\Vert\overline{u_{n}}%
\Vert_{L^{2}}\Vert\rho_{n+1}\Vert_{L^{\infty}}\Vert\overline{\rho_{n+1}}%
\Vert_{L^{2}}\\
\leq\,  &  CPA_{0}\left(  \Vert\overline{\omega_{n}}\Vert_{\dot{B}_{2,q}%
^{s-2}\cap\dot{H}^{-1}}+\Vert\overline{\omega_{n+1}}\Vert_{\dot{B}_{2,q}%
^{s-2}\cap\dot{H}^{-1}}\right)  \left(  \Vert\overline{\omega_{n+1}}%
\Vert_{\dot{H}^{-1}}^{2}+\Vert\overline{\rho_{n+1}}\Vert_{L^{2}}^{2}\right)
^{\frac{1}{2}}.
\end{split}
\]
Integrating over $(0,t)$ and using a Gr\"{o}nwall-type inequality, we have
that
\begin{equation}%
\begin{split}
\Vert\overline{\omega_{n+1}}(t)\Vert_{\dot{H}^{-1}}+\Vert\overline{\rho_{n+1}%
}(t)\Vert_{L^{2}}\leq\,  &  \Vert\overline{\omega_{n+1}}(0)\Vert_{\dot{H}%
^{-1}}+\Vert\overline{\rho_{n+1}}(0)\Vert_{L^{2}}\\
&  +CPA_{0}\int_{0}^{t}\Vert\overline{\omega_{n}}(\tau)\Vert_{\dot{B}%
_{2,q}^{s-2}\cap\dot{H}^{-1}}+\Vert\overline{\omega_{n+1}}(\tau)\Vert_{\dot
{B}_{2,q}^{s-2}\cap\dot{H}^{-1}}\,d\tau.
\end{split}
\label{H-1L2}%
\end{equation}
On the other hand, we applies $\Delta_{j}$ to the first and second equations
in (\ref{BCcurl2approx_dif}), and afterwards take the $\dot{H}^{-1}$-product
with $\Delta_{j}\overline{\omega_{n+1}}$ and the $L^{2}$-product with
$\Delta_{j}\overline{\rho_{n+1}}$ in order to get
\[%
\begin{split}
\langle\partial_{t}\Delta_{j}\omega_{n+1},\Delta_{j}\omega_{n+1}\rangle
_{\dot{H}^{-1}}=\,  &  -\langle\Delta_{j}(\overline{u_{n}}\cdot\nabla
)\omega_{n+1},\Delta_{j}\omega_{n+1}\rangle_{\dot{H}^{-1}}+\langle\lbrack
u_{n-1}\cdot\nabla,\Lambda^{-1}\Delta_{j}]\overline{\omega_{n+1}},\Lambda
^{-1}\Delta_{j}\omega_{n+1}\rangle_{L^{2}}\\
&  +\kappa\langle\Delta_{j}\partial_{1}\overline{\rho_{n+1}},\Delta_{j}%
\omega_{n+1}\rangle_{\dot{H}^{-1}},\\
\langle\partial_{t}\Delta_{j}\rho_{n+1},\Delta_{j}\rho_{n+1}\rangle_{L^{2}%
}=\,  &  -\langle\Delta_{j}(\overline{u_{n}}\cdot\nabla)\rho_{n+1},\Delta
_{j}\rho_{n+1}\rangle_{L^{2}}+\langle\lbrack u_{n-1}\cdot\nabla,\Delta
_{j}]\overline{\rho_{n+1}},\Delta_{j}\rho_{n+1}\rangle_{L^{2}}\\
&  +\kappa\langle\Delta_{j}\overline{u_{2,n+1}},\Delta_{j}\rho_{n+1}%
\rangle_{L^{2}}.
\end{split}
\]
Adding the two previous equalities, employing the property
\[
\langle\Delta_{j}\partial_{1}\overline{\rho_{n+1}},\Delta_{j}\omega
_{n+1}\rangle_{\dot{H}^{-1}}+\langle\Delta_{j}\overline{u_{2,n+1}},\Delta
_{j}\rho_{n+1}\rangle_{L^{2}}=0,
\]
and using Cauchy-Schwarz and H\"{o}lder inequalities, we obtain that
\[%
\begin{split}
\frac{1}{2}\frac{d}{dt}  &  \left(  \Vert\Delta_{j}\overline{\omega_{n+1}%
}\Vert_{\dot{H}^{-1}}^{2}+\Vert\Delta_{j}\overline{\rho_{n+1}}\Vert_{L^{2}%
}^{2}\right)  \leq\,\Bigl(\Vert\Delta_{j}(\overline{u_{n}}\cdot\nabla
)\omega_{n+1}\Vert_{\dot{H}^{-1}}+\Vert\lbrack u_{n-1}\cdot\nabla,\Lambda
^{-1}\Delta_{j}]\overline{\omega_{n+1}}\Vert_{L^{2}}\\
&  +\Vert\Delta_{j}(\overline{u_{n}}\cdot\nabla)\rho_{n+1}\Vert_{L^{2}}%
+\Vert\lbrack u_{n-1}\cdot\nabla,\Delta_{j}]\overline{\rho_{n+1}}\Vert_{L^{2}%
}\Bigr)\left(  \Vert\Delta_{j}\overline{\omega_{n+1}}\Vert_{\dot{H}^{-1}}%
^{2}+\Vert\Delta_{j}\overline{\rho_{n+1}}\Vert_{L^{2}}^{2}\right)  ^{\frac
{1}{2}}.
\end{split}
\]
Integrating over $(0,t)$ and using a Gr\"{o}nwall-type inequality, it follows
that
\[%
\begin{split}
\Vert\Delta_{j}\overline{\omega_{n+1}}(t)\Vert_{\dot{H}^{-1}}  &  +\Vert
\Delta_{j}\overline{\rho_{n+1}}(t)\Vert_{L^{2}}\leq\Vert\Delta_{j}%
\overline{\omega_{n+1}}(0)\Vert_{\dot{H}^{-1}}+\Vert\Delta_{j}\overline
{\rho_{n+1}}(0)\Vert_{L^{2}}\\
&  +\int_{0}^{t}\Vert\Delta_{j}(\overline{u_{n}}(\tau)\cdot\nabla)\omega
_{n+1}(\tau)\Vert_{\dot{H}^{-1}}\,d\tau+\int_{0}^{t}\Vert\lbrack u_{n-1}%
(\tau)\cdot\nabla,\Lambda^{-1}\Delta_{j}]\overline{\omega_{n+1}}(\tau
)\Vert_{L^{2}}\,d\tau\\
&  +\int_{0}^{t}\Vert\Delta_{j}(\overline{u_{n}}(\tau)\cdot\nabla)\rho
_{n+1}(\tau)\Vert_{L^{2}}\,d\tau+\int_{0}^{t}\Vert\lbrack u_{n-1}(\tau
)\cdot\nabla,\Delta_{j}]\overline{\rho_{n+1}}(\tau)\Vert_{L^{2}}\,d\tau.
\end{split}
\]
Taking into account Lemma \ref{inequalitybernstein}, multiplying by
$2^{(s-1)j}$ and taking the $l^{q}(\mathbb{Z})$-norm, we have that
\begin{equation}%
\begin{split}
&  \Vert\overline{\omega_{n+1}}(t)\Vert_{\dot{B}_{2,q}^{s-2}}+\Vert
\overline{\rho_{n+1}}(t)\Vert_{\dot{B}_{2,q}^{s-1}}\leq C\Vert\overline
{\omega_{n+1}}(0)\Vert_{\dot{B}_{2,q}^{s-2}}+C\Vert\overline{\rho_{n+1}%
}(0)\Vert_{\dot{B}_{2,q}^{s-1}}\\
&  +C\int_{0}^{t}\Vert(\overline{u_{n}}(\tau)\cdot\nabla)\omega_{n+1}%
(\tau)\Vert_{\dot{B}_{2,q}^{s-2}}\,d\tau+C\int_{0}^{t}\left(  \sum
_{j\in\mathbb{Z}}2^{(s-1)jq}\Vert\lbrack u_{n-1}(\tau)\cdot\nabla,\Lambda
^{-1}\Delta_{j}]\overline{\omega_{n+1}}(\tau)\Vert_{L^{2}}^{q}\right)
^{\frac{1}{q}}\,d\tau\\
&  +C\int_{0}^{t}\Vert(\overline{u_{n}}(\tau)\cdot\nabla)\rho_{n+1}(\tau
)\Vert_{\dot{B}_{2,q}^{s-1}}\,d\tau+C\int_{0}^{t}\left(  \sum_{j\in\mathbb{Z}%
}2^{(s-1)jq}\Vert\lbrack u_{n-1}(\tau)\cdot\nabla,\Delta_{j}]\overline
{\rho_{n+1}}(\tau)\Vert_{L^{2}}^{q}\right)  ^{\frac{1}{q}}\,d\tau\\
&  :=I_{1}+I_{2}+I_{3}+I_{4}+I_{5}+I_{6}.
\end{split}
\label{Bs-2Bs-1}%
\end{equation}
First, thanks to the embedding $\dot{B}_{2,q}^{s-1}\cap\dot{H}^{-1}%
\hookrightarrow L^{2}$ and the inequality
\begin{equation}
\Vert\overline{\omega_{n+1}}(0)\Vert_{\dot{B}_{2,q}^{s-2}\cap\dot{H}^{-1}%
}+\Vert\overline{\rho_{n+1}}(0)\Vert_{B_{2,q}^{s-1}}\leq C2^{-n}(\Vert
\omega_{0}\Vert_{\dot{B}_{2,q}^{s-1}\cap\dot{H}^{-1}}+\Vert\rho_{0}%
\Vert_{B_{2,q}^{s-1}})=CA_{0}2^{-n}, \label{I1I2}%
\end{equation}
we can estimate $I_{1}$ and $I_{2}$. For $I_{3}$ and $I_{5}$, let us first
note that by the equality $(f\cdot\nabla)g=\nabla\cdot(f\otimes g)$ for all
divergence free vector fields $f$, H\"{o}lder inequality, Remark
\ref{remark1}, the equality $u_{n}=\nabla^{\perp}(-\Delta)^{-1}\omega_{n}$,
the embedding $\dot{B}_{2,q}^{s-1}\cap\dot{H}^{-1}\hookrightarrow
B_{2,q}^{s-1}$ and estimate (\ref{unif_est}), we have that
\[%
\begin{split}
\Vert\nabla\cdot(\overline{u_{n}}\otimes\omega_{n+1})\Vert_{\dot{B}%
_{2,q}^{s-2}}  &  \leq C\Vert\overline{u_{n}}\otimes\omega_{n+1}\Vert_{\dot
{B}_{2,q}^{s-1}}\leq C(\Vert\overline{u_{n}}\Vert_{L^{\infty}}\Vert
\omega_{n+1}\Vert_{\dot{B}_{2,q}^{s-1}}+\Vert\omega_{n+1}\Vert_{L^{\infty}%
}\Vert\overline{u_{n}}\Vert_{\dot{B}_{2,q}^{s-1}})\\
&  \leq2CPA_{0}\Vert\overline{\omega_{n}}\Vert_{\dot{B}_{2,q}^{s-1}\cap\dot
{H}^{-1}},\\
\Vert(\overline{u_{n}}\cdot\nabla)\rho_{n+1}\Vert_{\dot{B}_{2,q}^{s-1}}  &
\leq C(\Vert\overline{u_{n}}\Vert_{L^{\infty}}\Vert\nabla\rho_{n+1}\Vert
_{\dot{B}_{2,q}^{s-1}}+\Vert\nabla\rho_{n+1}\Vert_{L^{\infty}}\Vert
\overline{u_{n}}\Vert_{\dot{B}_{2,q}^{s-1}})\\
&  \leq2CPA_{0}\Vert\overline{\omega_{n}}\Vert_{\dot{B}_{2,q}^{s-1}\cap\dot
{H}^{-1}}.
\end{split}
\]
Then,
\begin{equation}
I_{3}+I_{5}\leq4CPA_{0}\int_{0}^{t}\Vert\overline{\omega_{n}}(\tau)\Vert
_{\dot{B}_{2,q}^{s-1}\cap\dot{H}^{-1}}\,d\tau. \label{I3I5}%
\end{equation}
In view of Remark \ref{remark2}, Lemma \ref{2.4} and using the same arguments
to estimate $I_{3}$ and $I_{5}$, we see that
\[%
\begin{split}
I_{4}  &  \leq C\int_{0}^{t}\Vert\nabla u_{n-1}(\tau)\Vert_{L^{\infty}}%
\Vert\overline{\omega_{n+1}}(\tau)\Vert_{\dot{B}_{2,q}^{s-2}}+\Vert
\overline{\omega_{n+1}}(\tau)\Vert_{L^{2}}\Vert u_{n-1}(\tau)\Vert_{\dot
{B}_{\infty,q}^{s-1}}\,d\tau\\
&  \leq2CPA_{0}\int_{0}^{t}\Vert\overline{\omega_{n+1}}(\tau)\Vert_{\dot
{B}_{2,q}^{s-2}\cap\dot{H}^{-1}}\,d\tau,\\
I_{6}  &  \leq C\int_{0}^{t}\Vert\nabla u_{n-1}(\tau)\Vert_{L^{\infty}}%
\Vert\overline{\rho_{n+1}}(\tau)\Vert_{\dot{B}_{2,q}^{s-1}}+\Vert
\overline{\rho_{n+1}}(\tau)\Vert_{L^{\infty}}\Vert u_{n-1}(\tau)\Vert_{\dot
{B}_{2,q}^{s}}\,d\tau\\
&  \leq2CPA_{0}\int_{0}^{t}\Vert\overline{\rho_{n+1}}(\tau)\Vert
_{B_{2,q}^{s-1}}\,d\tau.
\end{split}
\]
Thus,
\begin{equation}
I_{4}+I_{6}\leq2CPA_{0}\int_{0}^{t}\Vert\overline{\omega_{n+1}}(\tau
)\Vert_{\dot{B}_{2,q}^{s-2}\cap\dot{H}^{-1}}+\Vert\overline{\rho_{n+1}}%
(\tau)\Vert_{B_{2,q}^{s-1}}\,d\tau. \label{I4I6}%
\end{equation}
Combining (\ref{H-1L2}) and (\ref{Bs-2Bs-1}), and using (\ref{I1I2}),
(\ref{I3I5}) and (\ref{I4I6}), it holds that
\begin{equation}%
\begin{split}
\Vert\overline{\omega_{n+1}}(t)\Vert_{\dot{B}_{2,q}^{s-2}\cap\dot{H}^{-1}%
}+\Vert\overline{\rho_{n+1}}(t)\Vert_{B_{2,q}^{s-1}}\leq\,  &  CA_{0}%
2^{-n}+5CPA_{0}\int_{0}^{t}\Vert\overline{\omega_{n}}(\tau)\Vert_{\dot
{B}_{2,q}^{s-2}\cap\dot{H}^{-1}}\,d\tau\\
&  +3CPA_{0}\int_{0}^{t}\Vert\overline{\omega_{n+1}}(\tau)\Vert_{\dot{B}%
_{2,q}^{s-2}\cap\dot{H}^{-1}}+\Vert\overline{\rho_{n+1}}(\tau)\Vert
_{B_{2,q}^{s-1}}\,d\tau.
\end{split}
\label{convergence}%
\end{equation}
By Gr\"{o}nwall inequality, it follows that
\[
\Vert\overline{\omega_{n+1}}(t)\Vert_{\dot{B}_{2,q}^{s-2}\cap\dot{H}^{-1}%
}+\Vert\overline{\rho_{n+1}}(t)\Vert_{B_{2,q}^{s-1}}\leq(CA_{0}2^{-n}%
+5CPA_{0}\int_{0}^{t}\Vert\overline{\omega_{n}}(\tau)\Vert_{\dot{B}%
_{2,q}^{s-2}\cap\dot{H}^{-1}}\,d\tau)e^{3CPA_{0}t}.
\]
Denoting $\overline{A_{n+1}}(t):=\Vert\overline{\omega_{n+1}}(t)\Vert_{\dot
{B}_{2,q}^{s-2}\cap\dot{H}^{-1}}+\Vert\overline{\rho_{n+1}}(t)\Vert
_{B_{2,q}^{s-1}}$, $f(t):=CA_{0}e^{3CPA_{0}t}$ and $g(t):=5CPA_{0}%
e^{3CPA_{0}t}$, we observe that
\[
\overline{A_{n+1}}(t)\leq2^{-n}f(t)+g(t)\int_{0}^{t}\overline{A_{n}}%
(\tau)\,d\tau.
\]
Following the iterative process, we arrive at
\[%
\begin{split}
\overline{A_{n+1}}(t)  &  \leq2^{-n}f(t)+g(t)\int_{0}^{t}\overline{A_{n}}%
(\tau)\,d\tau\\
&  \leq2^{-n}f(t)+g(t)\int_{0}^{t}\left(  2^{-(n-1)}f(\tau)+g(\tau)\int
_{0}^{\tau}\overline{A_{n-1}}(\tau^{\prime})\,d\tau^{\prime}\right)  \,d\tau\\
&  \leq2^{-n}f(t)+g(t)t\left(  2^{-(n-1)}f(t)+g(t)\int_{0}^{t}\overline
{A_{n-1}}(\tau)\,d\tau\right) \\
&  =2^{-n}f(t)\left(  1+2g(t)t\right)  +g(t)^{2}t\int_{0}^{t}\overline
{A_{n-1}}(\tau)\,d\tau\\
&  \leq2^{-n}f(t)\left(  1+2g(t)t+(2g(t)t)^{2}\right)  +g(t)^{3}t^{2}\int
_{0}^{t}\overline{A_{n-2}}(\tau)\,d\tau\\
&  \hspace{4cm}\vdots\\
&  \leq2^{-n}f(t)\sum_{i=0}^{n-1}(2g(t)t)^{i}+g(t)^{n}t^{n-1}\int_{0}%
^{t}\overline{A_{1}}(\tau)\,d\tau\\
&  \leq2^{-n}f(t)\sum_{i=0}^{n-1}(2g(t)t)^{i}+2PA_{0}g(t)^{n}t^{n}.
\end{split}
\]
Let $T^{\prime}\leq T$ be such that $g(t)t\leq\frac{1}{4}$ for all
$t\in(0,T^{\prime})$. Then, it is fulfilled that $\sum_{i=0}^{n-1}%
(2g(t)t)^{i}\leq\sum_{i=0}^{n-1}\frac{1}{2^{i}}<\infty$ and $(g(t)t)^{n}%
\leq4^{-n}$ for all $n\in\mathbb{N}$ and $t\in(0,T^{\prime})$. Moreover, there
exists $C_{3}>0$ such that $f(t)\leq C_{3}$ for all $t\in(0,T^{\prime})$.
Therefore, for some constant $C_{4}>0,$ we have
\[
\overline{A_{n+1}}(t)\leq C_{3}2^{-n}+2PA_{0}4^{-n}\leq C_{4}2^{-n}.
\]
Let $n,m\in\mathbb{N}$ be such that $n>m$, then there exists $C_{5}>0$
satisfying
\[
\Vert(\omega_{n}-\omega_{m})(t)\Vert_{\dot{B}_{2,q}^{s-2}\cap\dot{H}^{-1}%
}+\Vert(\rho_{n}-\rho_{m})(t)\Vert_{B_{2,q}^{s-1}}\leq\sum_{i=m}%
^{n-1}\overline{A_{i+1}}(t)\leq C_{5}\sum_{i=m}^{n-1}2^{-i}.
\]
This implies that $\left\{  \omega_{n}\right\}  _{n\in\mathbb{N}}$ and
$\left\{  \rho_{n}\right\}  _{n\in\mathbb{N}}$ are Cauchy in the spaces
$L^{\infty}(0,T^{\prime};\dot{B}_{2,q}^{s-2}\cap\dot{H}^{-1})$ and $L^{\infty
}(0,T^{\prime};B_{2,q}^{s-1})$, respectively. Therefore, there are $\omega\in
L^{\infty}(0,T^{\prime};\dot{B}_{2,q}^{s-2}\cap\dot{H}^{-1})$ and $\rho\in
L^{\infty}(0,T^{\prime};B_{2,q}^{s-1})$ such that if $n\rightarrow\infty$,
then
\[
\omega_{n}\longrightarrow\omega\,\text{ in }\,L^{\infty}(0,T^{\prime};\dot
{B}_{2,q}^{s-2}\cap\dot{H}^{-1})\,\text{ and }\,\rho_{n}\longrightarrow
\rho\,\text{ in }\,L^{\infty}(0,T^{\prime};B_{2,q}^{s-1}).
\]
Furthermore, since the sequences $\left\{  \omega_{n}\right\}  _{n\in
\mathbb{N}}$ and $\left\{  \rho_{n}\right\}  _{n\in\mathbb{N}}$ belong to
$C([0,T^{\prime}];\dot{B}_{2,q}^{s-1}\cap\dot{H}^{-1})$ and $C([0,T^{\prime
}];B_{2,q}^{s})$, respectively, we have that $\omega\in C([0,T^{\prime}%
];\dot{B}_{2,q}^{s-2}\cap\dot{H}^{-1})$ and $\rho\in C([0,T^{\prime}%
];B_{2,q}^{s-1})$.

On the other hand, as $\left\{  \omega_{n}\right\}  _{n\in\mathbb{N}}$ and
$\left\{  \rho_{n}\right\}  _{n\in\mathbb{N}}$ are bounded in $L^{\infty
}(0,T^{\prime};\dot{B}_{2,q}^{s-1}\cap\dot{H}^{-1})$ and $L^{\infty
}(0,T^{\prime};B_{2,q}^{s})$, respectively, we can extract subsequences
$\left\{  \omega_{n_{j}}\right\}  _{j\in\mathbb{N}}$ and $\left\{  \rho
_{n_{j}}\right\}  _{j\in\mathbb{N}}$ such that $\omega_{n_{j}}\overset{\ast
}{\rightharpoonup}\omega$ and $\rho_{n_{j}}\overset{\ast}{\rightharpoonup}%
\rho$ in $L^{\infty}(0,T^{\prime};\dot{B}_{2,q}^{s-1}\cap\dot{H}^{-1})$ and
$L^{\infty}(0,T^{\prime};B_{2,q}^{s})$, respectively. Thus,
\begin{equation}%
\begin{split}
\omega\in C([0,T^{\prime}];\dot{B}_{2,q}^{s-2}(\mathbb{R}^{2})\cap\dot{H}%
^{-1}(\mathbb{R}^{2}))  &  \cap L^{\infty}(0,T^{\prime};\dot{B}_{2,q}%
^{s-1}(\mathbb{R}^{2})\cap\dot{H}^{-1}(\mathbb{R}^{2})),\\
\rho\in C([0,T^{\prime}];B_{2,q}^{s-1}(\mathbb{R}^{2}))  &  \cap L^{\infty
}(0,T^{\prime};B_{2,q}^{s}(\mathbb{R}^{2})),
\end{split}
\label{solution_space0}%
\end{equation}
with $\Vert\omega\Vert_{L^{\infty}(0,T^{\prime};\dot{B}_{2,q}^{s-1}\cap\dot
{H}^{-1})},\Vert\rho\Vert_{L^{\infty}(0,T^{\prime};B_{2,q}^{s})}\leq PA_{0},$
where $P$ and $A_{0}$ are as in (\ref{unif_est}).

Now, with the above convergence in hand, we sketch the convergence of the
nonlinearity of the second equation and the coupling term in
(\ref{BCcurl2approx}). The others follow similarly and are left to the reader.
By H\"{o}lder's inequality, Remark \ref{remark1}, the identity $u_{n}%
-u=\nabla^{\perp}(-\Delta)^{-1}(\omega_{n}-\omega)$, the embedding
$B_{2,q}^{s-1}(\mathbb{R}^{2})\hookrightarrow B_{2,q}^{m+1}(\mathbb{R}^{2})$
and $B_{2,q}^{s-1}(\mathbb{R}^{2})\hookrightarrow B_{\infty,q}^{m}%
(\mathbb{R}^{2})$ with $0\leq m\leq s-2$, Lemma \ref{inequalitybernstein} and
(\ref{unif_est}), it follows that
\[%
\begin{split}
\int_{0}^{t}  &  \Vert(u_{n}(\tau)\cdot\nabla)\rho_{n+1}(\tau)-(u(\tau
)\cdot\nabla)\rho(\tau)\Vert_{B_{2,q}^{m+1}}\,d\tau\\
&  \leq\int_{0}^{t}\Vert((u_{n}-u)(\tau)\cdot\nabla)\rho_{n+1}(\tau
)\Vert_{B_{2,q}^{m+1}}+\Vert(u(\tau)\cdot\nabla)(\rho_{n+1}-\rho)(\tau
)\Vert_{B_{2,q}^{m+1}}\,d\tau\\
&  \leq CPA_{0}\int_{0}^{t}\Vert\omega_{n}-\omega\Vert_{B_{2,q}^{s-2}}%
+\Vert\rho_{n+1}-\rho\Vert_{B_{2,q}^{m+2}}+\Vert\rho_{n+1}-\rho\Vert
_{B_{2,2}^{1}}\,d\tau\\
&  \leq2CPA_{0}T(\Vert\omega_{n}-\omega\Vert_{L^{\infty}(0,T;B_{2,q}^{s-2}%
)}+\Vert\rho_{n+1}-\rho\Vert_{L^{\infty}(0,T;B_{2,q}^{m+2})}+\Vert\rho
_{n+1}-\rho\Vert_{L^{\infty}(0,T;B_{2,2}^{1})})\\
&  \rightarrow0\text{ as }n\rightarrow\infty,
\end{split}
\]
which implies
\[
\int_{0}^{t}(u_{n}(\tau)\cdot\nabla)\rho_{n+1}(\tau)\,d\tau\longrightarrow
\int_{0}^{t}(u(\tau)\cdot\nabla)\rho(\tau)\,d\tau\,\text{ in }\,B_{2,q}%
^{m+1}(\mathbb{R}^{2}),\,\text{ as }\,n\rightarrow\infty.
\]
Moreover, since $\rho_{n+1}\longrightarrow\rho$ in $L^{\infty}(0,T;B_{2,q}%
^{m+1}(\mathbb{R}^{2}))$, we have
\[
\left\Vert \int_{0}^{t}\kappa(\partial_{1}\rho_{n+1}(\tau)\,-\partial_{1}%
\rho(\tau))\,d\tau\right\Vert _{L^{\infty}(0,T;\dot{B}_{2,q}^{m}\cap\dot
{H}^{-1})}\leq C|\kappa|T\left\Vert \rho_{n+1}\,-\rho\right\Vert _{L^{\infty
}(0,T;B_{2,q}^{m+1})}\,\longrightarrow0,\text{ as }n\rightarrow\infty,\,\text{
}%
\]
and then
\[
\int_{0}^{t}\kappa\partial_{1}\rho_{n+1}(\tau)\,d\tau\longrightarrow\int
_{0}^{t}\kappa\partial_{1}\rho(\tau)\,d\tau\text{ in }\dot{B}_{2,q}%
^{m}(\mathbb{R}^{2})\cap\dot{H}^{-1}(\mathbb{R}^{2}).
\]
So, we can pass the limit in the integral formulation of approximate system
(\ref{BCcurl2approx}) and obtain that $(\omega,\rho)$ is an integral solution
for (\ref{Bouss1}) in $L^{\infty}(0,T;\dot{B}_{2,q}^{m}(\mathbb{R}^{2}%
)\cap\dot{H}^{-1}(\mathbb{R}^{2}))\times L^{\infty}(0,T;B_{2,q}^{m+1}%
(\mathbb{R}^{2}))$, namely
\begin{equation}%
\begin{split}
&  \omega(t)-\omega_{0}=\int_{0}^{t}(u(\tau)\cdot\nabla)\omega(\tau
)+\kappa\partial_{1}\rho(\tau)\,d\tau,\\
&  \rho(t)-\rho_{0}=\int_{0}^{t}(u(\tau)\cdot\nabla)\rho(\tau)+\kappa
u_{2}(\tau)\,d\tau.
\end{split}
\label{aux-int-sol-1}%
\end{equation}
Considering $y_{k}=\dot{S}_{k}\omega,$ $z_{k}=S_{k}\rho,$ using
(\ref{solution_space0}), and proceeding as in the proof of
(\ref{aux-cont-7000}) and (\ref{continuity_seq}), it follows that
$\{y_{k}\}_{k\in\mathbb{N}}\subset C([0,T];\dot{B}_{2,q}^{s-1}(\mathbb{R}%
^{2})\cap\dot{H}^{-1}(\mathbb{R}^{2}))$, $\{z_{k}\}_{k\in\mathbb{N}}\subset
C([0,T];B_{2,q}^{s}(\mathbb{R}^{2}))$, $y_{k}\rightarrow\omega$ in $L^{\infty
}(0,T;\dot{B}_{2,q}^{s-1}(\mathbb{R}^{2})\cap\dot{H}^{-1}(\mathbb{R}^{2}))$
and $z_{k}\rightarrow\rho$ in $L^{\infty}(0,T;B_{2,q}^{s}(\mathbb{R}^{2}))$,
as $k\rightarrow\infty$. Then, $\omega\in$ $C([0,T];\dot{B}_{2,q}%
^{s-1}(\mathbb{R}^{2})\cap\dot{H}^{-1}(\mathbb{R}^{2}))$, $\rho\in
C([0,T];B_{2,q}^{s}(\mathbb{R}^{2})),$ and the integral system
(\ref{aux-int-sol-1}) is indeed verified in $(\dot{B}_{2,q}^{s-2}%
(\mathbb{R}^{2})\cap\dot{H}^{-1}(\mathbb{R}^{2}))\times B_{2,q}^{s-1}%
(\mathbb{R}^{2})$, and consequently $\omega\in C^{1}([0,T];\dot{B}_{2,q}%
^{s-2}(\mathbb{R}^{2})\cap\dot{H}^{-1}(\mathbb{R}^{2}))$ and $\rho\in
C^{1}([0,T];B_{2,q}^{s-1}(\mathbb{R}^{2}))$, as desired.

\textbf{Uniqueness.} In this part, we suppose that system (\ref{Bouss1})
possesses two solutions $(\omega^{1},\rho^{1})$ and $(\omega^{2},\rho^{2})$
with the same initial data $(\omega_{0},\rho_{0})$ and we show that
$\omega^{1}=\omega^{2}$ and $\rho^{1}=\rho^{2}$. For that, we set
$\widetilde{\omega}:=\omega^{2}-\omega^{1}$ and $\widetilde{\rho}:=\rho
^{2}-\rho^{1}$, respectively. Then, $(\widetilde{\omega},\widetilde{\rho})$
satisfy the following system
\begin{equation}
\left\{
\begin{split}
&  \partial_{t}\widetilde{\omega}+(\widetilde{u}\cdot\nabla)\omega^{2}%
+(u^{1}\cdot\nabla)\widetilde{\omega}=\kappa\partial_{1}\widetilde{\rho},\\
&  \partial_{t}\widetilde{\rho}+(\widetilde{u}\cdot\nabla)\rho^{2}+(u^{1}%
\cdot\nabla)\widetilde{\rho}=\kappa\widetilde{u_{2}},\\
&  \widetilde{u}=\nabla^{\perp}(-\Delta)^{-1}\widetilde{\omega},\\
&  \widetilde{\omega}\mid_{t=0}=0,\,\widetilde{\rho}\mid_{t=0}=0.
\end{split}
\right.  \label{system_uniqueness}%
\end{equation}
Considering the spaces $L^{\infty}(0,T;\dot{B}_{2,q}^{s-2}\cap\dot{H}^{-1})$
and $L^{\infty}(0,T;B_{2,q}^{s-1})$ in (\ref{system_uniqueness}) and employing
the argument used in (\ref{BCcurl2approx_dif}) to estimate $\widetilde{\omega
}$ and $\widetilde{\rho}$, we obtain an inequality similar to inequality
(\ref{convergence}) as follows
\[
\Vert\widetilde{\omega}(t)\Vert_{\dot{B}_{2,q}^{s-2}\cap\dot{H}^{-1}}%
+\Vert\widetilde{\rho}(t)\Vert_{B_{2,q}^{s-1}}\leq CPA_{0}\int_{0}^{t}%
\Vert\widetilde{\omega}(\tau)\Vert_{\dot{B}_{2,q}^{s-2}\cap\dot{H}^{-1}}%
+\Vert\widetilde{\rho}(\tau)\Vert_{B_{2,q}^{s-1}}\,d\tau.
\]
Therefore, by Gr\"{o}nwall's inequality we obtain $\Vert\widetilde{\omega
}(t)\Vert_{\dot{B}_{2,q}^{s-2}\cap\dot{H}^{-1}}+\Vert\widetilde{\rho}%
(t)\Vert_{B_{2,q}^{s-1}}\leq0$ for all $t\in\lbrack0,T]$. Thus, $\Vert
\widetilde{\omega}\Vert_{L^{\infty}(0,T;\dot{B}_{2,q}^{s-2}\cap\dot{H}^{-1}%
)}=\Vert\widetilde{\rho}\Vert_{L^{\infty}(0,T;B_{2,q}^{s-1})}=0$, implying
that $\omega^{1}=\omega^{2}$ and $\rho^{1}=\rho^{2}$, which shows the
uniqueness of solution to (\ref{Bouss1}).

\section{Long-time solvability}

In this section we prove the long-time solvability of (\ref{Bouss1}) for large
values of $\left\vert \kappa\right\vert $. We start with a proposition
containing a blow-up criterion.

\begin{Proposition}
Let $s$ and $q$ be such that $s>2$ with $1\leq q\leq\infty$ or $s=2$ with
$q=1$. For $\omega_{0}\in\dot{B}_{2,q}^{s-1}(\mathbb{R}^{2})\cap\dot{H}%
^{-1}(\mathbb{R}^{2})$ and $\rho_{0}\in B_{2,q}^{s}(\mathbb{R}^{2})$, consider
$(\omega,\rho)$ the corresponding solution of (\ref{Bouss1}) satisfying
\[%
\begin{split}
\omega\in C([0,T);\dot{B}_{2,q}^{s-1}(\mathbb{R}^{2})\cap\dot{H}%
^{-1}(\mathbb{R}^{2}))  &  \cap C^{1}([0,T);\dot{B}_{2,q}^{s-2}(\mathbb{R}%
^{2})\cap\dot{H}^{-1}(\mathbb{R}^{2})),\\
\rho\in C([0,T);B_{2,q}^{s}(\mathbb{R}^{2}))  &  \cap C^{1}([0,T);B_{2,q}%
^{s-1}(\mathbb{R}^{2})),
\end{split}
\]
where $T>0$ is an existence time. If $\int_{0}^{T}\Vert\nabla\rho
(t)\Vert_{L^{\infty}}+\Vert\nabla u(t)\Vert_{L^{\infty}}\,dt<\infty$, then
there exists $T^{\prime}>T$ such that $(\omega,\rho)$ can be extended to
$[0,T^{\prime})$ with
\[%
\begin{split}
\omega\in C([0,T^{\prime});\dot{B}_{2,q}^{s-1}(\mathbb{R}^{2})\cap\dot{H}%
^{-1}(\mathbb{R}^{2}))  &  \cap C^{1}([0,T^{\prime});\dot{B}_{2,q}%
^{s-2}(\mathbb{R}^{2})\cap\dot{H}^{-1}(\mathbb{R}^{2})),\\
\rho\in C([0,T^{\prime});B_{2,q}^{s}(\mathbb{R}^{2}))  &  \cap C^{1}%
([0,T^{\prime});B_{2,q}^{s-1}(\mathbb{R}^{2})).
\end{split}
\]

\end{Proposition}

\textbf{Proof.} By standard procedures used to estimate $\omega$ and $\rho$ in
Besov norms, we obtain the following estimates:
\begin{equation}%
\begin{split}
\Vert &  \omega(t)\Vert_{\dot{B}_{2,q}^{s-1}\cap\dot{H}^{-1}}+\Vert
\rho(t)\Vert_{B_{2,q}^{s}}\leq C\Vert\omega(0)\Vert_{\dot{B}_{2,q}^{s-1}%
\cap\dot{H}^{-1}}+C\Vert\rho(0)\Vert_{B_{2,q}^{s}}\\
&  +\int_{0}^{t}\left(  \sum_{j\in\mathbb{Z}}2^{sjq}\Vert\lbrack u(\tau
)\cdot\nabla,\Lambda^{-1}\Delta_{j}]\omega(\tau)\Vert_{L^{2}}^{q}\right)
^{\frac{1}{q}}\,d\tau+\int_{0}^{t}\left(  \sum_{j\in\mathbb{Z}}2^{sjq}%
\Vert\lbrack u(\tau)\cdot\nabla,\Delta_{j}]\rho(\tau)\Vert_{L^{2}}^{q}\right)
^{\frac{1}{q}}\,d\tau.
\end{split}
\label{blowup_estimate}%
\end{equation}
If we denote by $I$ and $J$ the penultimate and last term of the previous
inequality, respectively, by Lemma \ref{2.4}, we have that
\[%
\begin{split}
I &  \leq C(\Vert\nabla u\Vert_{L^{\infty}}\Vert\omega\Vert_{\dot{B}%
_{2,q}^{s-1}}+\Vert\omega\Vert_{L^{\infty}}\Vert u\Vert_{\dot{B}_{2,q}^{s}%
}\leq C\Vert\omega\Vert_{\dot{B}_{2,q}^{s-1}\cap\dot{H}^{-1}}\left(
\Vert\nabla u\Vert_{L^{\infty}}+\Vert\omega\Vert_{L^{\infty}}\right)  ,\\
J &  \leq C(\Vert\nabla u\Vert_{L^{\infty}}\Vert\rho\Vert_{\dot{B}_{2,q}^{s}%
}+\Vert\nabla\rho\Vert_{L^{\infty}}\Vert u\Vert_{\dot{B}_{2,q}^{s}}\leq
C\left(  \Vert\omega\Vert_{\dot{B}_{2,q}^{s-1}\cap\dot{H}^{-1}}+\Vert\rho
\Vert_{B_{2,q}^{s}}\right)  \left(  \Vert\nabla u\Vert_{L^{\infty}}%
+\Vert\omega\Vert_{L^{\infty}}\right)  .
\end{split}
\]
Denoting $z_{s,q}(t):=\Vert\omega(t)\Vert_{\dot{B}_{2,q}^{s-1}\cap\dot{H}%
^{-1}}+\Vert\rho(t)\Vert_{B_{2,q}^{s}}$, employing the above inequalities in
(\ref{blowup_estimate}) and using $\left\Vert \omega\right\Vert _{L^{\infty}%
}\leq2\left\Vert \nabla u\right\Vert _{L^{\infty}}$, it holds
\[
z_{s,q}(t)\leq Cz_{s,q}(0)+C\int_{0}^{t}z_{s,q}(\tau)(\Vert\nabla\rho
(\tau)\Vert_{L^{\infty}}+\Vert\nabla u(\tau)\Vert_{L^{\infty}})\,d\tau.
\]
In turn, Gr\"{o}nwall-type inequality leads us to
\begin{equation}
z_{s,q}(t)\leq z_{s,q}(0)\exp{\left(  C_{6}\int_{0}^{t}\Vert\nabla\rho
(\tau)\Vert_{L^{\infty}}+\Vert\nabla u(\tau)\Vert_{L^{\infty}}d\tau\right)
},\label{inequality_blowupcriterion}%
\end{equation}
for all $t\in\lbrack0,T)$, where $C_{6}>0$ is a constant. Thus, by standard
arguments, $(\omega,\rho)$ can be continued to $[0,T^{\prime}]$ for some
$T^{\prime}>T$, whenever
\[
\int_{0}^{T}\Vert\nabla\rho(t)\Vert_{L^{\infty}}+\Vert\nabla u(t)\Vert
_{L^{\infty}}\,dt<\infty.
\]
\begin{flushright}$\blacksquare$\end{flushright}

\textbf{Long-time solvability. }For $\omega_{0}\in\dot{B}_{2,q}^{s}%
(\mathbb{R}^{2})\cap\dot{H}^{-1}(\mathbb{R}^{2})$, $\rho_{0}\in B_{2,q}%
^{s+1}(\mathbb{R}^{2})$, let $(\omega,\rho)$ be the solution of system
(\ref{Bouss1}) satisfying
\[%
\begin{split}
\omega\in C([0,T^{\ast});\dot{B}_{2,q}^{s}(\mathbb{R}^{2})\cap\dot{H}%
^{-1}(\mathbb{R}^{2})) &  \cap C^{1}([0,T^{\ast});\dot{B}_{2,q}^{s-1}%
(\mathbb{R}^{2})\cap\dot{H}^{-1}(\mathbb{R}^{2})),\\
\rho\in C([0,T^{\ast});B_{2,q}^{s+1}(\mathbb{R}^{2})) &  \cap C^{1}%
([0,T^{\ast});B_{2,q}^{s}(\mathbb{R}^{2})),
\end{split}
\]
with maximal existence time $T^{\ast}>0.$ If $T^{\ast}=\infty$, we are done.
Assume that $T^{\ast}<\infty.$ Denoting $V^{\pm}:=\omega\pm\Lambda\rho$, we
can use Duhamel's principle to get
\[
V^{\pm}(t)=e^{\pm\kappa\mathcal{R}_{1}t}V_{0}-\int_{0}^{t}e^{\pm
\kappa\mathcal{R}_{1}(\tau-t)}(f\pm\Lambda g)(\tau)\ d\tau,
\]
where $f=(u\cdot\nabla)\omega$ and $g=(u\cdot\nabla)\rho$. For $0\leq t\leq
T^{\ast},$ we define
\[
\mathcal{M}(t):=\int_{0}^{t}\Vert V^{\pm}(\tau)\Vert_{\dot{B}_{\infty,1}^{0}%
}\,d\tau.
\]
In what follows, we continue to use the notation $z_{s+1,q}(0)=\Vert\omega
_{0}\Vert_{\dot{B}_{2,q}^{s}\cap\dot{H}^{-1}}+\Vert\rho_{0}\Vert
_{B_{2,q}^{s+1}}$. We first consider the case $s=2$ with $q=1$. We can
estimate
\begin{align*}
\mathcal{M}(t) &  \leq C\int_{0}^{t}\Vert e^{\pm\kappa\mathcal{R}_{1}\tau
}V_{0}\Vert_{\dot{B}_{\infty,1}^{0}}\ d\tau+C\int_{0}^{t}\left\Vert \int
_{0}^{\tau}e^{\pm\kappa\mathcal{R}_{1}(\tau^{\prime}-\tau)}(f\pm\Lambda
g)(\tau^{\prime})\,d\tau^{\prime}\right\Vert _{\dot{B}_{\infty,1}^{0}}%
\ d\tau\\
&  :=K_{1}+K_{2}.
\end{align*}
For $K_{1}$, we use H\"{o}lder's inequality and Lemma
\ref{theorem-besov-strichartz} with $r=\infty$ to get
\[%
\begin{split}
K_{1} &  \leq Ct^{1-\frac{1}{\gamma}}\Vert e^{\pm\kappa\mathcal{R}_{1}\cdot
}V_{0}\Vert_{L^{\gamma}(0,\infty;\dot{B}_{\infty,1}^{0})}\\
&  \leq Ct^{1-\frac{1}{\gamma}}\,|\kappa|^{-\frac{1}{\gamma}}\Vert V_{0}%
\Vert_{\dot{B}_{2,1}^{1}}\\
&  \leq Ct^{1-\frac{1}{\gamma}}\,|\kappa|^{-\frac{1}{\gamma}}z_{2,1}(0).
\end{split}
\]
For $K_{2}$, we employ the Minkowski and H\"{o}lder inequalities, Lemma
\ref{theorem-besov-strichartz} and Remark \ref{remark1} to obtain
\[%
\begin{split}
K_{2} &  \leq C\int_{0}^{t}\int_{0}^{\tau}\Vert e^{\pm\kappa\mathcal{R}%
_{1}(\tau^{\prime}-\tau)}(f\pm\Lambda g)(\tau^{\prime})\Vert_{\dot{B}%
_{\infty,1}^{0}}\,d\tau^{\prime}\,d\tau\\
&  =C\int_{0}^{t}\int_{\tau^{\prime}}^{t}\Vert e^{\pm\kappa\mathcal{R}%
_{1}(\tau^{\prime}-\tau)}(f\pm\Lambda g)(\tau^{\prime})\Vert_{\dot{B}%
_{\infty,1}^{0}}\,d\tau\,d\tau^{\prime}\\
&  \leq Ct^{1-\frac{1}{\gamma}}\int_{0}^{t}\Vert e^{\pm\kappa\mathcal{R}%
_{1}(\tau^{\prime}-\cdot)}(f\pm\Lambda g)(\tau^{\prime})\Vert_{L^{\gamma}%
(\tau^{\prime},t;\dot{B}_{\infty,1}^{0})}\,d\tau^{\prime}\\
&  \leq Ct^{1-\frac{1}{\gamma}}\,|\kappa|^{-\frac{1}{\gamma}}\int_{0}^{t}%
\Vert(f\pm\Lambda g)(\tau^{\prime})\Vert_{\dot{B}_{2,1}^{1}}\,d\tau^{\prime}\\
&  \leq Ct^{1-\frac{1}{\gamma}}\,|\kappa|^{-\frac{1}{\gamma}}\int_{0}%
^{t}(\Vert\omega(\tau)\Vert_{\dot{B}_{2,1}^{0}}\Vert\omega(\tau)\Vert_{\dot
{B}_{2,1}^{2}}+\Vert\omega(\tau)\Vert_{\dot{B}_{2,1}^{0}}\Vert\rho(\tau
)\Vert_{\dot{B}_{2,1}^{3}}+\Vert\rho(\tau)\Vert_{\dot{B}_{2,1}^{2}}\Vert
\omega(\tau)\Vert_{\dot{B}_{2,1}^{1}})\,d\tau\\
&  \leq Ct^{1-\frac{1}{\gamma}}\,|\kappa|^{-\frac{1}{\gamma}}\int_{0}%
^{t}(\Vert\omega(\tau)\Vert_{\dot{B}_{2,1}^{2}\cap\dot{H}^{-1}}^{2}+\Vert
\rho(\tau)\Vert_{B_{2,1}^{3}}^{2})\,d\tau.
\end{split}
\]
Thus, for each $0<t<T^{\ast}$, we use (\ref{inequality_blowupcriterion}), the
embedding $\dot{B}_{\infty,1}^{0}\hookrightarrow L^{\infty}$ and the equality
$u=\nabla^{\perp}(-\Delta)^{-1}\omega$ to get
\[%
\begin{split}
\mathcal{M}(t) &  \leq Ct^{1-\frac{1}{\gamma}}|\kappa|^{-\frac{1}{\gamma}%
}\left(  z_{2,1}(0)+\int_{0}^{t}\Vert\omega(\tau)\Vert_{\dot{B}_{2,1}^{2}%
\cap\dot{H}^{-1}}^{2}+\Vert\rho(\tau)\Vert_{B_{2,1}^{3}}^{2}\,d\tau\right)  \\
&  \leq Ct^{1-\frac{1}{\gamma}}\,|\kappa|^{-\frac{1}{\gamma}}\left(
z_{3,1}(0)+z_{3,1}(0)^{2}\int_{0}^{t}e^{C_{6}\mathcal{M}(\tau)}d\tau\right)
\\
&  \leq Ct^{1-\frac{1}{\gamma}}|\kappa|^{-\frac{1}{\gamma}}z_{3,1}(0)\left(
1+z_{3,1}(0)te^{C_{6}\mathcal{M}(t)}\right)  .
\end{split}
\]
Now, we deal with the case $s>2$ with $1\leq q\leq\infty$. For each $1\leq
q\leq\infty,$ we take $2<r\leq\infty$ such that $q\leq r$. Note that since
$s-1>1$ we have the nonhomogeneous embedding $B_{\infty,\infty}^{s-1}%
\hookrightarrow W^{1,\infty}$, so that we can estimate $\Vert V\Vert_{\dot
{B}_{\infty,1}^{0}}$ by the $B_{\infty,\infty}^{s-2}$-norm of $V$. Thus, we
have that%
\[%
\begin{split}
\mathcal{M}(t) &  \leq C\int_{0}^{t}\Vert e^{\kappa\mathcal{R}_{1}\tau}%
V_{0}\Vert_{B_{\infty,\infty}^{s-2}}\ d\tau+C\int_{0}^{t}\left\Vert \int
_{0}^{\tau}e^{\kappa\mathcal{R}_{1}(\tau^{\prime}-\tau)}(f+\Lambda
g)(\tau^{\prime})\,d\tau^{\prime}\right\Vert _{B_{\infty,\infty}^{s-2}}%
\ d\tau\\
&  :=K_{3}+K_{4}.
\end{split}
\]
We now estimate $K_{3}$. Using the embedding $\dot{B}_{\infty,q}%
^{s-2}\hookrightarrow\dot{B}_{\infty,\infty}^{s-2}$, H\"{o}lder's inequality
and Lemma \ref{theorem-besov-strichartz}, we can estimate%
\[%
\begin{split}
\int_{0}^{t}\Vert e^{\kappa\mathcal{R}_{1}\tau}V_{0}\Vert_{\dot{B}%
_{\infty,\infty}^{s-2}}\ d\tau &  \leq\int_{0}^{t}\Vert e^{\kappa
\mathcal{R}_{1}\tau}V_{0}\Vert_{\dot{B}_{\infty,q}^{s-2}(\mathbb{R}^{2}%
)}\,d\tau\\
&  \leq t^{1-\frac{1}{\gamma}}\Vert e^{\kappa\mathcal{R}_{1}\tau}V_{0}%
\Vert_{L^{\gamma}(0,\infty;\dot{B}_{\infty,q}^{s-2})}\\
&  \leq Ct^{1-\frac{1}{\gamma}}\,|\kappa|^{-\frac{1}{\gamma}}\Vert V_{0}%
\Vert_{\dot{B}_{2,q}^{s-1}}\\
&  \leq Ct^{1-\frac{1}{\gamma}}\,|\kappa|^{-\frac{1}{\gamma}}z_{s,q}(0).
\end{split}
\]
Also, by Lemma \ref{LemmaStrichartz}, we have that%
\[%
\begin{split}
\int_{0}^{t}\Vert e^{\kappa\mathcal{R}_{1}\tau}V_{0}\Vert_{L^{\infty}}\ d\tau
&  \leq t^{1-\frac{1}{\gamma}}\Vert e^{\kappa\mathcal{R}_{1}\tau}V_{0}%
\Vert_{L^{\gamma}(0,\infty;L^{\infty})}\\
&  \leq Ct^{1-\frac{1}{\gamma}}|\kappa|^{-\frac{1}{\gamma}}\Vert V_{0}%
\Vert_{L^{2}}.
\end{split}
\]
Therefore,
\[
K_{3}\leq Ct^{1-\frac{1}{\gamma}}\,|\kappa|^{-\frac{1}{\gamma}}z_{s,q}(0).
\]
We proceed to estimate $K_{4}$. H\"{o}lder's inequality, the embedding
$\dot{B}_{\infty,q}^{s-2}\hookrightarrow\dot{B}_{\infty,\infty}^{s-2}$ and
Lemma \ref{theorem-besov-strichartz} yield
\[%
\begin{split}
\int_{0}^{t}\left\Vert \int_{0}^{\tau}e^{\kappa\mathcal{R}_{1}(\tau^{\prime
}-\tau)}(f+\Lambda g)(\tau^{\prime})\,d\tau^{\prime}\right\Vert _{\dot
{B}_{\infty,\infty}^{s-2}}\ d\tau &  \leq C\int_{0}^{t}\int_{0}^{\tau}\Vert
e^{\kappa\mathcal{R}_{1}(\tau^{\prime}-\tau)}(f+\Lambda g)(\tau^{\prime}%
)\Vert_{\dot{B}_{\infty,q}^{s-2}}\,d\tau^{\prime}\,d\tau\\
&  =C\int_{0}^{t}\int_{\tau^{\prime}}^{t}\Vert e^{\kappa\mathcal{R}_{1}%
(\tau^{\prime}-\tau)}(f+\Lambda g)(\tau^{\prime})\Vert_{\dot{B}_{\infty
,q}^{s-2}}\,d\tau\,d\tau^{\prime}\\
&  \leq Ct^{1-\frac{1}{\gamma}}\int_{0}^{t}\Vert e^{\kappa\mathcal{R}_{1}%
(\tau^{\prime}-\tau)}(f+\Lambda g)(\tau^{\prime})\Vert_{L^{\gamma}%
(\tau^{\prime},t;\dot{B}_{\infty,q}^{s-2})}\,d\tau^{\prime}\\
&  \leq Ct^{1-\frac{1}{\gamma}}\,|\kappa|^{-\frac{1}{\gamma}}\int_{0}^{t}%
\Vert(f+\Lambda g)(\tau)\Vert_{\dot{B}_{2,q}^{s-1}}\,d\tau\\
&  \leq Ct^{1-\frac{1}{\gamma}}\,|\kappa|^{-\frac{1}{\gamma}}\int_{0}%
^{t}(\Vert\omega(\tau)\Vert_{\dot{B}_{2,q}^{s}\cap\dot{H}^{-1}}^{2}+\Vert
\rho(\tau)\Vert_{B_{2,q}^{s+1}}^{2})\,d\tau.
\end{split}
\]
Here we have proceeded as for the term $K_{2}$ in the case $s=2$, $q=1$. Also,
by Lemma \ref{LemmaStrichartz}, the continuity of $\mathcal{R}_{l}$ in $L^{2}$
and the embedding $B_{2,q}^{s+1}\hookrightarrow B_{2,q}^{s}\hookrightarrow
L^{2}$, we arrive at
\[%
\begin{split}
\int_{0}^{t}\left\Vert \int_{0}^{\tau}e^{\kappa\mathcal{R}_{1}(\tau^{\prime
}-\tau)}(f+\Lambda g)(\tau^{\prime})\,d\tau^{\prime}\right\Vert _{L^{\infty}%
}\,d\tau^{\prime} &  \leq Ct^{1-\frac{1}{\gamma}}|\kappa|^{-\frac{1}{\gamma}%
}\int_{0}^{t}\Vert(f+\Lambda g)(\tau^{\prime})\Vert_{L^{2}}\,d\tau^{\prime}\\
&  \leq Ct^{1-\frac{1}{\gamma}}|\kappa|^{-\frac{1}{\gamma}}\int_{0}^{t}%
(\Vert\omega(\tau)\Vert_{\dot{B}_{2,q}^{s}\cap\dot{H}^{-1}}^{2}+\Vert\rho
(\tau)\Vert_{B_{2,q}^{s+1}}^{2})\,d\tau.
\end{split}
\]
It follows that
\[
K_{4}\leq C\,t^{1-\frac{1}{\gamma}}|\kappa|^{-\frac{1}{\gamma}}\int_{0}%
^{t}\Vert\omega(\tau)\Vert_{\dot{B}_{2,q}^{s}\cap\dot{H}^{-1}}^{2}+\Vert
\rho(\tau)\Vert_{B_{2,q}^{s+1}}^{2}\,d\tau^{\prime}.
\]
Thus, for each $0<t<T^{\ast}$, we have
\[%
\begin{split}
\mathcal{M}(t) &  \leq C\,t^{1-\frac{1}{\gamma}}|\kappa|^{-\frac{1}{\gamma}%
}\left(  z_{s,q}(0)+\int_{0}^{t}\Vert\omega(\tau)\Vert_{\dot{B}_{2,q}^{s}%
\cap\dot{H}^{-1}}^{2}+\Vert\rho(\tau)\Vert_{B_{2,q}^{s+1}}^{2}\,d\tau\right)
\\
&  \leq Ct^{1-\frac{1}{\gamma}}|\kappa|^{-\frac{1}{\gamma}}\left(
z_{s+1,q}(0)+z_{s+1,q}(0)^{2}\int_{0}^{t}e^{C_{6}\mathcal{M}(\tau)}%
d\tau\right)  \\
&  \leq Ct^{1-\frac{1}{\gamma}}|\kappa|^{-\frac{1}{\gamma}}z_{s+1,q}(0)\left(
1+z_{s+1,q}(0)te^{C_{6}\mathcal{M}(t)}\right)  .
\end{split}
\]
Therefore, for both cases of $s$ and $q$ such that $s=2$ with $q=1$ or $s>2$
with $1\leq q\leq\infty$, we have that there exists $C_{7}>0$ such that
\begin{equation}
\mathcal{M}(t)\leq C_{7}t^{1-\frac{1}{\gamma}}|\kappa|^{-\frac{1}{\gamma}%
}z_{s+1,q}(0)\left(  1+z_{s+1,q}(0)te^{C_{6}\mathcal{M}(t)}\right)
.\label{aux-time-ext-10}%
\end{equation}
Next, for each $0<T<\infty$ we define $\tilde{T}=\sup D_{T}$, where
\[
D_{T}=\{t\in\lbrack0,T]\cap\lbrack0,T^{\ast})\ |\ \mathcal{M}(t)\leq
C_{7}T^{1-\frac{1}{\gamma}}z_{s+1,q}(0)\}.
\]
We first show that $\tilde{T}=\min\{T,T^{\ast}\}.$ We proceed by
contradiction. So, assume on the contrary that $\tilde{T}<\min\{T,T^{\ast}\}$.
We have that there exists $T_{1}$ such that $\tilde{T}<T_{1}<\min\{T,T^{\ast
}\}$. It follows that
\[%
\begin{split}
\omega\in C([0,T_{1}];\dot{B}_{2,q}^{s-1}(\mathbb{R}^{2})\cap\dot{H}%
^{-1}(\mathbb{R}^{2})) &  \cap C^{1}([0,T_{1}];\dot{B}_{2,q}^{s-2}%
(\mathbb{R}^{2})\cap\dot{H}^{-1}(\mathbb{R}^{2})),\\
\rho\in C([0,T_{1}];B_{2,q}^{s+1}(\mathbb{R}^{2})) &  \cap C^{1}%
([0,T_{1}];B_{2,q}^{s}(\mathbb{R}^{2})),
\end{split}
\]
$\mathcal{M}(t)$ is uniformly continuous on $[0,T_{1}],$ and
\begin{equation}
\mathcal{M}(\tilde{T})\leq C_{7}\,T^{1-\frac{1}{\gamma}}z_{s+1,q}%
(0).\label{time ext 2}%
\end{equation}
We now take a $|\kappa|$ large enough so that
\begin{equation}
|\kappa|^{\frac{1}{\gamma}}\geq2\left(  1+(\Vert\omega_{0}\Vert_{\dot{B}%
_{2,q}^{s}\cap\dot{H}^{-1}}+\Vert\rho_{0}\Vert_{B_{2,q}^{s+1}})T\exp{\left(
C_{6}C_{7}T^{1-\frac{1}{\gamma}}(\Vert\omega_{0}\Vert_{\dot{B}_{2,q}^{s}%
\cap\dot{H}^{-1}}+\Vert\rho_{0}\Vert_{B_{2,q}^{s+1}})\right)  }\right)
.\label{time ext 3}%
\end{equation}
Using (\ref{aux-time-ext-10}), (\ref{time ext 2}) and (\ref{time ext 3}), we
obtain that
\[%
\begin{split}
\mathcal{M}(\tilde{T}) &  \leq C_{7}(\tilde{T})^{1-\frac{1}{\gamma}}%
|\kappa|^{-\frac{1}{\gamma}}z_{s+1,q}(0)\left(  1+z_{s+1,q}(0)\tilde{T}%
\exp\left(  C_{6}\mathcal{M}(\tilde{T})\right)  \right)  \\
&  \leq C_{7}T^{1-\frac{1}{\gamma}}z_{s+1,q}(0)|\kappa|^{-\frac{1}{\gamma}%
}\left(  1+z_{s+1,q}(0)T\exp\left(  C_{6}C_{7}T^{1-\frac{1}{\gamma}}%
z_{s+1,q}(0)\right)  \right)  \\
&  \leq\frac{1}{2}C_{7}\,T^{1-\frac{1}{\gamma}}z_{s+1,q}(0).
\end{split}
\]
Thus, we can choose $T_{2}$ such that $\tilde{T}<T_{2}<T_{1}$ with
$\mathcal{M}(T_{2})\leq C_{7}T^{1-\frac{1}{\gamma}}z_{s,q}(0)$. This
contradicts the definition of $\tilde{T}$. It follows that $\tilde{T}%
=\min\{T,T^{\ast}\}$ when $\kappa$ verifies (\ref{time ext 3}). If $T^{\ast
}<T$, then $T^{\ast}=\tilde{T}=\sup D_{T}$ and
\[
\mathcal{M}(t)\leq C_{7}T^{1-\frac{1}{\gamma}}z_{s+1,q}(0)<\infty,\text{ for
all }0\leq t<T^{\ast}.
\]
It follows that $\int_{0}^{T^{\ast}}\Vert\nabla\rho(t)\Vert_{L^{\infty}}%
+\Vert\nabla u(t)\Vert_{L^{\infty}}\,dt\leq C\mathcal{M}(T^{\ast})<\infty$
and, in view of the blow-up criterion, we obtain a contradiction with the
maximality of $T^{\ast}$. This concludes the proof.\begin{flushright}$\blacksquare$\end{flushright}


\begin{thebibliography}{99}                                                                                               %


\bibitem {AngFerr}V. Angulo-Castillo and L. C. F. Ferreira, On the 3D Euler
equations with Coriolis force in borderline Besov spaces. Commun. Math. Sci 16
(1) (2018), 145--164.

\bibitem {AngFerrKos}V. Angulo-Castillo, L. C. F. Ferreira, and L. Kosloff,
Long-time solvability for the 2D dispersive SQG equation with improved
regularity. Discrete Contin. Dyn. Syst. 40 (3) (2020), 1411--1433.

\bibitem {BaMaNi2}A. Babin, A. Mahalov, and B. Nicolaenko, 3D Navier-Stokes
and Euler equations with initial data characterized by uniformly large
vorticity. Indiana Univ. Math. J. 50 (2001), Special Issue, 1--35.

\bibitem {Bahouri-Chemin94}H. Bahouri and J.-Y. Chemin, \'{E}quations de
transport relatives \'{a} des champs de vecteurs non-lipschitziens et
m\'{e}canique des fluides. Arch. Rational Mech. Anal. 127 (2) (1994), 159--181.

\bibitem {BahouriCheminDanchin}H. Bahouri, J. Chemin, and R. Danchin, Fourier
Analysis and Nonlinear Partial Differential Equations, Grundlehren der
mathematischen Wissenschaften, vol. 343, Springer, Heidelberg (2011).

\bibitem {Barbu}V. Barbu, Differential equations, Translated from the
Romanian, Originally published by Junimea, Springer, Ia{\c{s}}i, 1985.

\bibitem {Wu Stab}O. Ben Said, U. Pandey and J. Wu, The stabilizing effect of
the temperature on buoyancy-driven fluids, Indiana Univ. Math. J. 71 (2022), 2605--2645.

\bibitem {Berg-Lofstrom}J. Bergh and J. L\"{o}fstr\"{o}m, Interpolation
Spaces. An Introduction, Springer-Verlag, Berlin-New York, 1976.

\bibitem {Bourgain-Li}J. Bourgain and D. Li, Strong ill-posedness of the
incompressible Euler equation in borderline Sobolev spaces. Invent. Math.
201(1) (2014), 97--157.

%\bibitem {Cannone}M. Cannone, C. Miao, and L. Xue, Global regularity for the
%supercritical dissipative quasi-geostrophic equation with large dispersive
%forcing. Proc. Lond. Math. Soc. 106 (2013), 650--674.


\bibitem {CCL}A. Castro, D. C\'{o}rdoba, and D. Lear, On the asymptotic
stability of stratified solutions for the 2D Boussinesq equations with a
velocity damping term. Math. Models Methods Appl. Sci. 29 (7) (2019), 1227--1277.

\bibitem {Chae2004}D. Chae, Local existence and blow-up criterion for the
Euler equations in the Besov spaces. Asymptot. Anal. 38 (2004), 339--358.

\bibitem {Chae}D. Chae, Global regularity for the 2D Boussinesq equations with
partial viscosity terms, Adv. in Math., 203, 497-513, 2006. .

%\bibitem {ChKN}D. Chae, S. Kim, and H. Nam, Local existence and blow-up
%criterion of H\"{o}lder continuous solutions of the Boussinesq equations.
%Nagoya Math. J. 155 (1999), 55--80.


\bibitem {Chemin92}J.-Y. Chemin, R\'{e}gularit\'{e} de la trajectoire des
particules d'un fluide parfait incompressible remplissant l'espace. J. Math.
Pures Appl. 71 (5) (1992), 407--417.

%\bibitem {CD99}P. Constantin and C.R. Doering, Infinite Prandtl number
%convection. J. Statistical Physics 94 (1999), 159--172.


%\bibitem {DoG}C. R. Doering and J. D. Gibbon, Applied Analysis of the
%Navier-Stokes Equations. Cambridge University Press, 1995.


\bibitem {D-P}R. Danchin and M. Paicu, Le th\'{e}or\`{e}me de Leray et le
th\'{e}or\`{e}me de Fujita-Kato pour le syst\`{e}me de Boussinesq
partiellement visqueux, Bull. Soc. Math. France, 136, 261-309, 2008.

\bibitem {D-W-Z-Z}Charles R. Doering, Jiahong Wu, Kun Zhao, and Xiaoming
Zheng, Long time behavior of the two-dimensional Boussinesq equations without
buoyancy diffusion. Phys.D, 376/377:144-159, 2018.

\bibitem {Dutrifoy}A. Dutrifoy, Examples of dispersive effects in non-viscous
rotating fluids. J. Math. Pures Appl. 84 (2005), 331--356.

\bibitem {EJJ}T. M. Elgindi and I.-J. Jeong, Finite-time singularity formation
for strong solutions to the Boussinesq system. Ann. PDE 6 (1) (2020), Paper
No. 5, 50 pp.

\bibitem {ElgMasm}T. M. Elgindi and N. Masmoudi, $L^{\infty}$ Ill-Posedness
for a Class of Equations Arising in Hydrodynamics. Arch Rational Mech Anal 235
(2020), 1979--2025.

\bibitem {EW}T. M. Elgindi and K. Widmayer, Sharp decay estimates for an
anisotropic linear semigroup and applications to the surface quasi-geostrophic
and inviscid Boussinesq systems. SIAM J. Math. Anal. 47 (2015), 4672--4684.

\bibitem {LCF-EVR}L.C.F. Ferreira and E.J. Villamizar-Roa, Strong solutions
and inviscid limit for Boussinesq system with partial viscosity. Commun. Math.
Sci. 11 (2) (2013), 421--439.

\bibitem {Gill}A.E. Gill, Atmosphere-Ocean Dynamics. International Geophysics
Series, vol. 30. Academic Press, Cambridge 1982.
%\bibitem {Fujii-2021}M. Fujii, Long time existence and asymptotic behavior of
%solutions for the 2D quasi-geostrophic equation with large dispersive forcing.
%J. Math. Fluid Mech. 23 (1) (2021), Paper No. 12, 19 pp.


%\bibitem {HH}Z. Hassainia and T. Hmidi, On the inviscid Boussinesq system with
%rough initial data. J. Math. Anal. Appl. 430 (2) (2015), 777--809.


\bibitem {H-K}T. Hmidi and S. Keraani, On the global well-possedness of the
two-dimensional Boussinesq system with a zero diffusivity, Adv. Diff. Eqs.,
12, 461-480, 2007

\bibitem {KLT}Y. Koh, S. Lee, and R. Takada, Strichartz estimates for the
Euler equations in the rotational framework. J. Differential Equations 256 (2)
(2014), 707--744.

\bibitem {Liu et al}X. Liu, M. Wang, and Z. Zhang, Local well-posedness and
blowup criterion of the Boussinesq equations in critical Besov spaces. J.
Math. Fluid Mech. 12 (2010), 280--292.

\bibitem {Majda}A. J. Majda, Introduction to PDEs and Waves for the Atmosphere
and Ocean. Courant Lecture Notes in Mathematics 9, AMS/CIMS, 2003.

\bibitem {PakPark}H. Pak and Y. Park, Existence of solution for the Euler
equations in a critical Besov space $B_{\infty,1}^{1}(R^{n})$. Comm. Partial
Differential Equations 29 (7-8) (2004), 1149--1166.

\bibitem {Pedlosky}J. Pedlosky, Geophysical Fluid Dynamics, Springer-Verlag,
New York, 1987.

\bibitem {Salmon}R. Salmon. Lectures on geophysical fluid dynamics. Oxford
University Press, New York, 1998.

\bibitem {SW}A. Stefanov and J. Wu, A global regularity result for the 2D
Boussinesq equations with critical dissipation. J. Anal. Math. 137 (1) (2019), 269--290.

\bibitem {T-W-Z-Z}L. Tao, J. Wu, K. Zhao, and X. Zheng, Stability near
hydrostatic equilibrium to the 2D Boussinesq equations without thermal
diffusion. Archive for Rational Mechanics and Analysis 237 (2) (2020), 585--630.

\bibitem {Takada2008}R. Takada, Local existence and blow-up criterion for the
Euler equations in Besov spaces of weak type. J. Evol. Equ. 8 (2008), 693--725.

\bibitem {Takada2021}R. Takada, Long time solutions for the 2D inviscid
Boussinesq equations with strong stratification. Manuscripta Math. 164 (1-2)
(2021), 223--250.

\bibitem {Vishik99}M. Vishik, Incompressible flows of an ideal fluid with
vorticity in borderline spaces of Besov type. Ann. Sci. \'{E}cole Norm. Sup.
32 (6) (1999), 769--812.

\bibitem {Wan-CCM-2020}R. Wan, Long time stability for the dispersive SQG
equation and Boussinesq equations in Sobolev space $H^{s}$. Commun. Contemp.
Math. 22 (3) (2020), 1850063, 13 pp.

\bibitem {WC17}R. Wan and J. Chen, Global well-posedness for the 2D dispersive
SQG equation and inviscid Boussinesq equations. Z. Angew. Math. Phys. 67 (4)
(2016), 22pp.

\bibitem {Wid}K. Widmayer, Convergence to stratified flow for an inviscid 3D
Boussinesq system. Commun. Math. Sci. 16 (6) (2018), 1713--1728.

\bibitem {WuRef}J. Wu, The 2d incompressible Boussinesq equations, Summer
School Lecture Notes, Peking University, 2012.

\bibitem {WuXuYe2015}J. Wu, X. Xu, and Z. Ye, Global smooth solutions to the
$n$-dimensional damped models of incompressible fluid mechanics with small
initial datum. J. Nonlinear Sci. 25 (1) (2015), 157--192.
\end{thebibliography}
\end{document}